\documentclass[a4paper,12pt]{amsart}
\usepackage{latexsym}
\usepackage{amsmath}
\usepackage{graphicx}
\usepackage{amscd}
\usepackage{amssymb}
\usepackage{times}
\usepackage{leftidx}

\newtheorem{thm}{Theorem}[section]
\newtheorem{prop}[thm]{Proposition}
\newtheorem{lem}[thm]{Lemma}
\newtheorem{cor}[thm]{Corollary}
\theoremstyle{definition}
\newtheorem{definition}[thm]{Definition}
\newtheorem{example}[thm]{Example}
\newtheorem{remark}[thm]{Remark}

\title[The Magnus representation and higher-order Alexander 
invariants]{The Magnus representation and higher-order Alexander 
invariants for homology cobordisms of surfaces} 
\author{Takuya SAKASAI}
\address{Graduate School of Mathematical Sciences, 
the University of Tokyo, 3-8-1 Komaba, Meguro-ku, 
Tokyo 153-8914, Japan}
\email{sakasai@ms.u-tokyo.ac.jp}


\newcommand{\Ker}{\mathop{\mathrm{Ker}}\nolimits}
\newcommand{\Hom}{\mathop{\mathrm{Hom}}\nolimits}
\newcommand{\Id}{\mathop{\mathrm{id}}\nolimits}
\newcommand{\Z}{\ensuremath{\mathbb{Z}}}
\newcommand{\Tmatrix}[1]{\mathop{\left( {#1} \right)}\nolimits}
\newcommand{\Aut}{\mathop{\mathrm{Aut}}\nolimits}
\newcommand{\IAut}{\mathop{\mathrm{IAut}}\nolimits}

\renewcommand{\Im}{\mathop{\mathrm{Im}}\nolimits}
\newcommand{\Mg}{\ensuremath{\mathcal{M}_{g,1}}}
\newcommand{\Cg}{\ensuremath{\mathcal{C}_{g,1}}}
\newcommand{\Hg}{\ensuremath{\mathcal{H}_{g,1}}}
\newcommand{\Lg}{\mathcal{L}_g}
\newcommand{\Slg}{\ensuremath{\mathcal{S}_{g}}}
\newcommand{\Sg}{\ensuremath{\Sigma_{g,1}}}
\newcommand{\Acy}{F^{\mathrm{acy}}}
\newcommand{\KK}{\mathcal{K}}
\newcommand{\ab}{\mathrm{ab}}
\newcommand{\sgn}{\mathop{\mathrm{sgn}}\nolimits}
\newcommand{\rank}{\mathop{\mathrm{rank}}\nolimits}
\newcommand{\zegamma}{\textstyle
\frac{\partial\zeta}{\partial \overrightarrow{\gamma}}}

\begin{document}

\begin{abstract}
The set of homology cobordisms from a surface to itself with markings 
of their boundaries has a 
natural monoid structure. 
To investigate the structure of this monoid, 
we define and study its Magnus representation and 
Reidemeister torsion invariants 
by generalizing Kirk-Livingston-Wang's argument over 
the Gassner representation of string links. 
Moreover, by applying Cochran and Harvey's 
framework of higher-order (non-commutative) Alexander invariants 
to them, we extract several pieces of information about 
the monoid and related objects. 
\end{abstract}

\maketitle

\section{Introduction}\label{sec:intro}
Let $\Sg$ be a compact connected 
oriented surface of genus $g \ge 1$ with one 
boundary component. 
A {\it homology cylinder} ({\it over} $\Sg$) consists of 
a homology cobordism from $\Sg$ to itself with 
markings of its boundary. We denote by $\Cg$ the set 
of isomorphisms classes of 
homology cylinders. Stacking two homology cylinders 
gives a new one, and by this, 
we can endow $\Cg$ with a monoid structure 
(see Section \ref{sec:def} for the 
precise definition). 
The origin of homology cylinders goes back to 
Habiro \cite{ha}, Garoufalidis-Levine \cite{gl} 
and Levine \cite{le}, 
where the clasper (or clover) surgery theory is 
effectively used to investigate the structure of $\Cg$. 

By a standard method, we can assign a homology cylinder 
to each homology 3-sphere or pure string link. Also, 
for a given homology cylinder, 
we can use an element of the mapping class group of $\Sg$
to construct another one by changing its markings. 
Since these operations preserve each monoid structure, 
$\Cg$ can be regarded as 
a simultaneous generalization of 
the monoid of homology 3-spheres, that of string links 
and the mapping class group, 
any of which plays an important role in the theory of 3-manifolds. 
On the other hand, there exists a natural way 
(called {\it closing}) 
to construct a closed 3-manifold from each homology cylinder. 
Therefore, through its monoid structure, $\Cg$ serves as 
an effective tool for classifying closed 3-manifolds. 

The aim of this paper is to study the structure of $\Cg$ 
from rather an algebraic point of view. We mainly use 
non-commutative rings arising from group rings to 
define some invariants such as 
the Magnus representation for $\Cg$ 
and Reidemeister torsion invariants. 
Note that our Magnus representation extends 
that for the mapping class group defined 
by Morita \cite{mo}, as 
the Gassner representation for string links due to 
Le Dimet \cite{ld} and  Kirk-Livingston-Wang \cite{klw} 
does that for the pure braid group. 
See Birman's book \cite{bi} for generalities of 
the ordinary (pre-extended) Magnus representation, 
including free differentials. 

After defining invariants using non-commutative rings, 
we shall need 
some devices to extract information from them. 
For that, we use the framework of higher-order Alexander 
invariants due to Cochran \cite{coc} and Harvey \cite{har,har2}. 
Higher-order Alexander invariants are 
those for finitely presentable groups 
interpreted as degrees of ``non-commutative Alexander polynomials'',  
which have some unclear ambiguity except their degrees. 
Historically, they are first defined for knot groups by Cochran, 
and then generalized for arbitrary finitely presentable groups 
by Harvey. 
Using them, Cochran and Harvey obtained various sharper 
results than those brought by 
the ordinary Alexander invariants --- lower bounds on the knot genus or 
the Thurston norm, necessary conditions for realizing a 
given group as the fundamental group of some compact oriented 
3-manifold, and so on. 
In the process of applying higher-order Alexander 
invariants to our case, we shall give its slight generalization 
(called {\it torsion-degree functions}) because of 
the difference of localizations of non-commutative rings used in 
the Magnus representation and higher-order Alexander invariants. 
Then we use it to study several properties of our invariants 
and relationships between them, from which we will obtain some 
information about the structure of $\Cg$ and related 3-manifolds. 

The outline of this paper is as follows. 
In Section \ref{sec:def}, we review the definition of homology cylinders 
as well as setting up our notation and terminology. 
Sections \ref{sec:magnus} and \ref{sec:gassner}, which are the first 
main part of this paper, 
are devoted to define the Magnus representation 
and study its fundamentals, including some examples. 
In Section \ref{sec:harvey}, we review the theory of higher-order 
Alexander invariants, following Harvey's papers \cite{har,har2}, 
and then define its generalization. 
In Section \ref{sec:apply}, which is the second main part, 
we observe several applications of our invariants.

\vspace{10pt}

\section{Homology cylinders}\label{sec:def}
\begin{figure}[htbp]
\begin{center}
\includegraphics{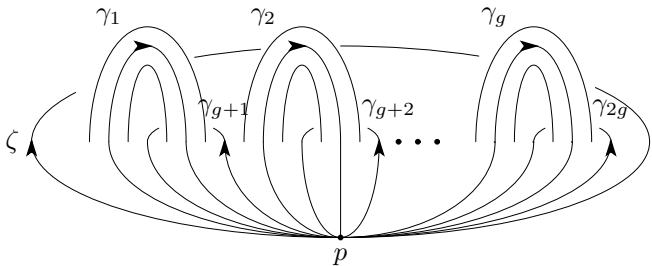}
\end{center}
\caption{A standard cell decomposition of $\Sigma_{g,1}$}
\label{fig:generator}
\end{figure}

Throughout the paper, we work in PL or smooth category. 
Let $\Sg$ be a compact connected oriented surface 
of genus $g \ge 1$ with one boundary component. We take a base point 
$p$ on the boundary of $\Sg$, and take $2g$ loops 
$\gamma_1, \ldots , \gamma_{2g}$ 
of $\Sg$ as shown in Figure \ref{fig:generator}. 
We consider them to be an embedded bouquet $R_{2g}$ of 
$2g$-circles tied at the base point $p \in \partial \Sg$. Then 
$R_{2g}$ and the boundary loop $\zeta$ of $\Sg$ 
together with one 2-cell make up a standard cell decomposition of $\Sg$. 
The fundamental group $\pi_1 \Sg$ of $\Sg$ is 
isomorphic to the free group $F_{2g}$ of rank $2g$ generated by 
$\gamma_1,\ldots,\gamma_{2g}$, in which 
$\zeta=\prod_{i=1}^g [\gamma_i,\gamma_{g+i}]$. 

\noindent
A {\it homology cylinder} $(M,i_+,i_-)$ ({\it over} $\Sg$), 
which has its origin in 
Habiro \cite{ha}, Garoufalidis-Levine \cite{gl} and 
Levine \cite{le}, consists of a 
compact oriented 3-manifold $M$ and two embeddings 
$i_+ , i_-: \Sg \rightarrow \partial M$ satisfying that 
\begin{enumerate}
\item $i_+$ is 
orientation-preserving and $i_-$ is orientation-reversing, 
\item $\partial M= i_+ (\Sg) \cup i_- (\Sg)$ and 
$i_+ (\Sg) \cap i_- (\Sg)=i_+ (\partial \Sg)
=i_- (\partial \Sg)$, 
\item $i_+ \bigl|_{\partial \Sg}=i_- \bigl|_{\partial \Sg}$, 
\item $i_+,i_-:H_\ast (\Sg) \rightarrow H_\ast (M)$ 
are isomorphisms. 
\end{enumerate}
\noindent
We denote $i_+ (p) =i_- (p)$ by $p \in \partial M$ again and consider 
it to be the base point of $M$. 
We write a homology cylinder by $(M,i_+,i_-)$ or simply by $M$. 

Two homology cylinders are said to be {\it isomorphic} if there exists 
an orientation-preserving diffeomorphism between 
the underlying 3-manifolds which is 
compatible with the embeddings of $\Sg$. 
We denote the set of isomorphism classes of homology cylinders 
by $\Cg$. Given two homology cylinders $M=(M,i_+,i_-)$ and 
$N=(N,j_+,j_-)$, we can construct a new 
homology cylinder $M \cdot N$ by 
\[M \cdot N = (M \cup_{i_- \circ (j_+)^{-1}} N, i_+,j_-).\]
Then $\Cg$ becomes a monoid with the unit 
$1_{\Cg}:=(\Sg \times I, \Id \times 1, \Id \times 0)$. 

From the monoid $\Cg$, we can construct the 
{\it homology cobordism group $\Hg$ 
of homology cylinders} as in the following way. 
Two homology cylinders $M=(M,i_+,i_-)$ and 
$N=(N,j_+,j_-)$ are {\it homology cobordant} if there exists 
a compact oriented 4-manifold $W$ such that 
\begin{enumerate}
\item $\partial W = M \cup (-N) /(i_+ (x)= j_+(x) , \,
i_- (x)=j_-(x)) \quad x \in \Sg$, 
\item the inclusions $M \hookrightarrow W$, $N \hookrightarrow W$ 
induce isomorphisms on the homology,  
\end{enumerate}
\noindent
where $-N$ is $N$ with opposite orientation. 
We denote by $\Hg$ 
the quotient set of $\Cg$ with respect to the equivalence relation 
of homology cobordism. The monoid structure of $\Cg$ induces 
a group structure of $\Hg$. In the group $\Hg$, the inverse of 
$(M,i_+,i_-)$ is given by $(-M,i_-,i_+)$.
\begin{example}\label{mgtocg}
For each element $\varphi$ of the mapping class group $\Mg$ of $\Sg$, 
we can construct a homology cylinder 
$M_{\varphi} \in \Cg$ by setting 
\[M_{\varphi}:=(\Sg \times I, \Id \times 1, 
\varphi \times 0),\]
where collars of $i_+ (\Sg)$ and $i_- (\Sg)$ are stretched half-way 
along $\partial \Sg \times I$. 
This gives injective homomorphisms 
$\Mg \hookrightarrow \Cg$ and $\Mg \hookrightarrow \Hg$. 
From this, we can regard $\Cg$ and $\Hg$ as enlargements of $\Mg$. 
\end{example}

Let $N_k (G):=G/(\Gamma^k G)$ be the $k$-th nilpotent 
quotient of a group $G$, where we define $\Gamma^1 G=G$ and 
$\Gamma^l G =[\Gamma^{l-1} G,G]$ for $l \ge 2$. 
For simplicity, we write 
$N_k (X)$ for $N_k (\pi_1 X)$ where 
$X$ is a connected topological space, 
and write $N_k$ for $N_k (F_{2g})=N_k (\Sg)$. 

Let $(M,i_+,i_-)$ be a homology cylinder. By definition, 
$i_+, i_- :\pi_1 \Sg \to \pi_1 M$ are both 2-{\it connected}, 
namely they induce isomorphisms on the first homology groups and 
epimorphisms on the second homology groups. 
Then, by Stallings' theorem \cite{st}, 
$i_+, i_- :N_k \xrightarrow{\cong} N_k (M)$ are isomorphisms 
for each $k \ge 2$. Using them, 
we obtain a monoid homomorphism 
\[\sigma_k: \Cg \longrightarrow \Aut N_k 
\qquad ((M,i_+,i_-) \mapsto (i_+)^{-1} \circ i_-).\]
It can be easily checked that $\sigma_k$ induces a group homomorphism 
$\sigma_k: \Hg \to \Aut N_k$. 
We define filtrations of $\Cg$ and $\Hg$ by 
\[\begin{array}{ll}
\Cg[1]:=\Cg, & \Cg[k]:=
\Ker \left( \Cg \xrightarrow{\sigma_k} \Aut N_k \right) \ 
\mbox{for $k \ge 2$}, \\
\Hg[1]:=\Hg, & \Hg[k]:=
\Ker \left( \Hg \xrightarrow{\sigma_k} \Aut N_k \right) \ 
\mbox{for $k \ge 2$}.
\end{array}\]

\vspace{10pt}

\section{The Magnus representation for homology cylinders}\label{sec:magnus}
We first summarize our notation. 
For a matrix $A$ with entries in a ring $R$, 
and a ring homomorphism $\varphi:R \to R'$, 
we denote by ${}^{\varphi} A$ the 
matrix obtained from $A$ by applying $\varphi$ to each entry. 
$A^T$ denotes the transpose of $A$. 
When $R=\Z G$ for a group $G$ or its right field of fractions 
(if exists), we denote by $\overline{A}$ 
the matrix obtained from $A$ by applying the involution induced 
from $(x \mapsto x^{-1},\ x \in G)$ to each entry. 

For a module $M$, we write $M^n$ and $M_n$ for the modules of 
column and row vectors with $n$ entries respectively. 

For a finite cell complex $X$ and 
its regular covering $X_{\Gamma}$ with respect to a homomorphism 
$\pi_1 X \to \Gamma$, $\Gamma$ acts on $X_{\Gamma}$ from 
the right through its deck transformation group. Therefore 
we regard the $\Z \Gamma$-cellular chain complex 
$C_{\ast} (X_{\Gamma})$ of $X_{\Gamma}$ as a collection of 
free right $\Z \Gamma$-modules consisting of column vectors together 
with differentials given by left multiplications of matrices. 
For each $\Z \Gamma$-bimodule $A$, 
the twisted chain complex $C_{\ast} (X;A)$ is given 
by the tensor product of 
the right $\Z \Gamma$-module $C_{\ast} (X_{\Gamma})$ and 
the left $\Z \Gamma$-module $A$, 
so that $C_{\ast} (X;A)$ and 
$H_{\ast} (X;A)$ are right $\Z \Gamma$-modules.

\subsection{Definition of the Magnus representation 
for homology cylinders}
In what follows, we fix an integer $k \ge 2$, which corresponds to 
the class of the nilpotent quotient. 
The following construction is based on Kirk-Livingston-Wang's work 
of the Gassner representation for string links in \cite{klw}. 

Let $(M,i_+,i_-) \in \Cg$ be a homology cylinder. 
By Stallings' theorem, 
$N_k$ and $N_k (M)$ are isomorphic. 
Since $N_k$ is a finitely generated torsion-free nilpotent group 
for each $k \ge 2$, we can embed $\Z N_k$ into the 
right field of fractions $\KK_{N_k} := \Z N_k (\Z N_k -\{ 0\})^{-1}$. 
(See Section \ref{sec:harvey}.) 
Similarly, we have $\Z N_k (M) \hookrightarrow 
\KK_{N_k (M)} := \Z N_k (M)(\Z N_k (M)-\{ 0\})^{-1}$. 
We consider the fields $\mathcal{K}_{N_k}$ and 
$\mathcal{K}_{N_k (M)}$ to be 
local coefficient systems on $\Sigma_{g,1}$ and $M$ respectively. 

By a standard argument using covering spaces 
(see for instance \cite[Proposition 2.1]{klw}, \cite[Lemma 5.11]{sa2}), 
we have the following.
\begin{lem}\label{relative}
$i_{\pm}: 
H_\ast (\Sg,p ;i_{\pm}^\ast \KK_{N_k (M)}) \to 
H_\ast (M,p ;\KK_{N_k (M)})$ are isomorphisms as 
right $\KK_{N_k (M)}$-vector spaces. 
\end{lem}
\begin{remark}\label{acyclic}
The same conclusion as in Lemma \ref{relative} can be drawn 
for the homology with coefficients in any 
$\Z \pi_1 (M)$-algebra $A$ satisfying the following property: 
Every matrix with entries in $\Z \pi_1 (M)$ sent to 
an invertible one by the augmentation map $\Z \pi_1 (M) \to \Z$ is 
also invertible in $A$. Note that 
$\KK_{N_k (M)}$ satisfies this property. 
\end{remark}
\noindent
Since $R_{2g} \subset \Sg$ is a deformation 
retract, we have
\[H_1 (\Sg,p;i_{\pm}^\ast \KK_{N_k (M)}) \cong 
H_1 (R_{2g},p;i_{\pm}^\ast \KK_{N_k (M)}) = 
C_1 (\widetilde{R_{2g}}) \otimes_{\pi_1 R_{2g}} 
i_{\pm}^\ast \KK_{N_k (M)} 
\cong \KK_{N_k (M)}^{2g}\]
with a basis
\[\{ \widetilde{\gamma_1} \otimes 1, \ldots , 
\widetilde{\gamma_{2g}} \otimes 1\} 
\subset C_1 (\widetilde{R_{2g}}) \otimes_{\pi_1 R_{2g}} 
i_{\pm}^\ast \KK_{N_k (M)}\] 
as a right $\KK_{N_k (M)}$-vector space, where 
$\widetilde{\gamma_i}$ is a lift 
of $\gamma_i$ on the universal covering $\widetilde{R_{2g}}$. 

\begin{definition}\label{def:Mag2}
$(1)$ For each $M=(M,i_+,i_-) \in \Cg$, we 
denote by $r'_k(M) \in GL(2g,\KK_{N_k (M)})$ 
the representation matrix of 
the right $\KK_{N_k (M)}$-isomorphism
\[\KK_{N_k (M)}^{2g} \cong 
H_1 (\Sg,p;i_-^\ast \KK_{N_k (M)}) 
\xrightarrow[i_-]{\cong} 
H_1 (M,p;\KK_{N_k (M)}) 
\xrightarrow[i_+^{-1}]{\cong} 
H_1 (\Sg,p;i_+^\ast \KK_{N_k (M)}) 
\cong \KK_{N_k (M)}^{2g}\]
$(2)$ The {\it Magnus representation} for $\Cg$ 
is the map $r_k :\Cg \to GL(2g,\KK_{N_k})$ 
which assigns to $M=(M,i_+,i_-) \in \Cg$ the 
matrix ${}^{i_+^{-1}} r'_k(M)$.
\end{definition}
\noindent
While we call $r_k(M)$ the Magnus 
``representation'', it is actually a crossed homomorphism, 
namely we have the following. 
\begin{thm} For $M_1, M_2 \in \Cg$, we have
\[r_k(M_1 \cdot M_2) = r_k(M_1) \cdot {}^{\sigma_k(M_1)} r_k(M_2).\]
\end{thm}
\begin{proof}
We write $M=M_1 \cdot M_2$ for simplicity. 
Let $i :M_1 \to M$ and $j :M_2 \to M$ 
be the natural inclusions. 
Since $M=(M,i \circ i_+, j \circ j_-)$ and 
$i \circ i_-=j \circ j_+$, the map 
\[H_1 (\Sg ,p;j_-^{\ast}j^{\ast}\KK_{N_k (M)}) 
\xrightarrow{j \circ j_-}
H_1 (M,p;\KK_{N_k (M)})
\xrightarrow{(i \circ i_+)^{-1}}
H_1 (\Sg ,p;i_+^{\ast}i^{\ast}\KK_{N_k (M)})\]
is given as the composite of 
\[H_1 (\Sg ,p;j_-^{\ast}j^{\ast}\KK_{N_k (M)}) 
\xrightarrow{j_-}
H_1 (M_2,p;j^{\ast}\KK_{N_k (M)})
\xrightarrow{j_+^{-1}}
H_1 (\Sg ,p;j_+^{\ast}j^{\ast}\KK_{N_k (M)})\]
and 
\[H_1 (\Sg ,p;i_-^{\ast}i^{\ast}\KK_{N_k (M)}) 
\xrightarrow{i_-}
H_1 (M_1,p;i^{\ast}\KK_{N_k (M)})
\xrightarrow{i_+^{-1}}
H_1 (\Sg ,p;i_+^{\ast}i^{\ast}\KK_{N_k (M)}).\]
Hence 
\[\begin{array}{crcl}
&r'_k(M) &=& 
{}^{i} r'_k(M_1) \cdot {}^{j} r'_k(M_2) \\
\Longrightarrow & {}^{(i \, i_+)^{-1}}r'_k(M) &=&
{}^{(i \, i_+)^{-1} i}r'_k(M_1) \cdot 
{}^{(i \, i_+)^{-1} j}r'_k(M_2) \\
 & &=&
{}^{i_+^{-1}}r'_k(M_1) \cdot
{}^{i_+^{-1} i^{-1} j}r'_k(M_2) \\
 & &=&
{}^{i_+^{-1}}r'_k(M_1) \cdot 
{}^{i_+^{-1} i_- j_+^{-1}}r'_k(M_2) \\
\Longrightarrow & r_k(M) &=&
r_k(M_1) \cdot {}^{\sigma_k(M_1)} r_k(M_2) 
\end{array}\]
This completes the proof. 
\end{proof}
\begin{thm}
$r_k:\Cg \rightarrow 
GL(2g,\KK_{N_k})$ factors through $\Hg$.
\end{thm}
\begin{proof} 
Suppose $M_1 =(M_1 ,i_+,i_-)$ and 
$M_2=(M_2 ,j_+,j_-) \in \Cg$ are homology cobordant by 
a homology cobordism $W$. Let 
$i:M_1 \to W$, $j:M_2 \to W$ be the natural inclusions. 
We may assume that $M_1 \cup M_2$ is a subcomplex of $W$ 
and that $W$ has only one 0-cell $p$. 
Since $\Z N_k \cong \Z N_k (W)$ 
by Stallings' theorem, we have 
$\KK_{N_k (W)} := \Z N_k (W)(\Z N_k (W)-\{ 0\})^{-1}$. 
We write $I_+:=i \circ i_+ = j \circ j_+$ and 
$I_-:=i \circ i_- = j \circ j_-$. 
Then we have the following commutative diagram:
\[\begin{CD}
H_1(\Sg,p;i_-^{\ast} i^{\ast}\KK_{N_k (W)}) 
@= H_1(\Sg,p;I_-^{\ast}\KK_{N_k (W)}) 
@= H_1(\Sg,p;j_-^{\ast}j^{\ast}\KK_{N_k (W)}) \\
@V{i_-}V{\cong}V @V{I_-}V{\cong}V @V{j_-}V{\cong}V \\
H_1(M_1 ,p;i^{\ast}\KK_{N_k (W)}) 
@>i>{\cong}> H_1(W,p;\KK_{N_k (W)}) 
@<j<{\cong}< H_1(M_2 ,p;j^{\ast}\KK_{N_k (W)})\\
@V{(i_+)^{-1}}V{\cong}V @V{(I_+)^{-1}}V{\cong}V 
@V{(j_+)^{-1}}V{\cong}V \\
H_1(\Sg,p;i_+^{\ast} i^{\ast}\KK_{N_k (W)}) 
@= H_1(\Sg,p;I_+^{\ast}\KK_{N_k (W)}) 
@= H_1(\Sg,p;j_+^{\ast} j^{\ast}\KK_{N_k (W)}) 
\end{CD}\]
The left vertical maps give 
${}^i r'_k (M_1)$ and the right ones 
give ${}^j r'_k (M_2)$. Applying $I_+^{-1}$, 
we obtain $r_k(M_1) = r_k(M_2)$. 
\end{proof}

\noindent
Consequently, we obtain the Magnus representation 
$r_k: \Hg \to GL(2g,\KK_{N_k})$, 
which is a crossed homomorphism. 
If we restrict $r_k$ 
to $\Cg[k]$ (resp.~ $\Hg[k]$), it becomes a monoid (resp.~group) 
homomorphism.
\begin{example}\label{fox}
For $\varphi \in \Mg \hookrightarrow \Aut F_{2g}$, we can obtain 
\[r_k(M_{\varphi}) = 
\overline{
\sideset{^{\rho_k}\!}{}
{\Tmatrix{
\displaystyle\frac{\partial \varphi(\gamma_j)}{\partial \gamma_i}
}}
}_{i,j},\]
where $\rho_k:\Z F_{2g} \to \Z N_k \subset \KK_{N_k}$ 
is the natural map and $\partial/\partial \gamma_i$ are 
free differentials. From this, we see that 
$r_k$ generalizes the original Magnus representation for 
$\Mg$ in \cite{mo}.
\end{example}

\subsection{Computation of the Magnus matrix}\label{subsec:computation}
In \cite{klw}, the Gassner matrix of a string link is computed from 
the Wirtinger presentation of the fundamental group of 
its exterior, which gives a 
finite presentation whose deficiency coincides with the number 
of strings. Recall that the deficiency of a finite presentation 
$P=\{x_1,\ldots,x_n \mid r_1,\ldots,r_m \} $ of a finitely 
presentable group $G$ is $n-m$, and the deficiency of $G$ is 
the maximum of all over the deficiencies of finite presentations of $G$. 
In our context, we do not have such a useful method in general. 
\begin{definition}\label{admissible}
For $(M,i_+,i_-) \in \Cg$, 
a presentation of $\pi_1 M$ is said to be {\it admissible} if 
it is of the form
\[\langle i_- (\gamma_1),\ldots,i_- (\gamma_{2g}), 
z_1 ,\ldots, z_l, 
i_+ (\gamma_1),\ldots,i_+ (\gamma_{2g}) \mid 
r_1, \ldots, r_{2g+l}
\rangle.\]
\end{definition}
Note that there does exist an admissible presentation 
for each homology cylinder $(M,i_+,i_-)$. 
Indeed, take a Morse function with no critical points of 
indices 0 and 3. Then $M$ can be seen as $\Sg \times I$ with 
some 1- and 2-handles. Since the Euler characteristics of 
$\Sg \times I$ and $M$ are the same, the numbers of the attached 
1- and 2- handles are the same. 
Therefore the presentation of $\pi_1 M$ obtained from 
a presentation of $\pi_1 (\Sg \times I) = F_{2g}$ 
with deficiency $2g$ by adding 
new generators and relations corresponding to the 
1- and 2-handles has deficiency $2g$ again. 
Our claim follows from this. 
(See also Section \ref{subsec:torsion}.)

Given an admissible presentation of $\pi_1 M$ 
as in Definition \ref{admissible}, 
we define $2g \times (2g+l)$, $l \times (2g+l)$ and 
$2g \times (2g+l)$ matrices $A,B,C$ by 
\[A=\overline{
\left(\frac{\partial r_j}{\partial i_-(\gamma_i)}
\right)}_{\begin{subarray}{c}
{}1 \le i \le 2g\\
1 \le j \le 2g+l
\end{subarray}}, \quad 
B=\overline{
\left(\frac{\partial r_j}{\partial z_i}
\right)}_{\begin{subarray}{c}
{}1 \le i \le l\\
1 \le j \le 2g+l
\end{subarray}}, \quad 
C=\overline{
\left(\frac{\partial r_j}{\partial i_+(\gamma_i)}
\right)}_{\begin{subarray}{c}
{}1 \le i \le 2g\\
1 \le j \le 2g+l
\end{subarray}}\]
at $\Z N_k (M)$, namely we apply the natural map 
\[\Z \langle i_- (\gamma_1),\ldots,i_- (\gamma_{2g}), 
z_1 ,\ldots, z_l, 
i_+ (\gamma_1),\ldots,i_+ (\gamma_{2g}) \rangle 
\to \Z \pi_1 (M) \to \Z N_k (M)\]
to each entry of the matrices obtained by free differentials.
\begin{prop}\label{howto}
$(1)$ The square matrix 
$\begin{pmatrix} A \\ B \end{pmatrix}$ 
is invertible as a matrix with entries 
in $\KK_{N_k (M)}$.\\
$(2)$ As matrices with entries 
in $\KK_{N_k (M)}$, we have 
\[(r'_k (M) \quad Z) \begin{pmatrix} A \\ B \end{pmatrix}
=-C\]
for some $2g \times l$ matrix $Z$.
\end{prop}
\begin{proof}
(1) Let $\mathfrak{t}:\Z N_k (M) \to \Z$ 
be the augmentation map. 
$\sideset{^\mathfrak{t}}{}
{\Tmatrix{
\begin{smallmatrix} A \\ B \end{smallmatrix}}}\! \cdot$ 
gives 
a presentation matrix of $H_1 (M) / \Phi_+$, where $\Phi_+$ 
is the subgroup of $H_1 (M)$ generated 
by $i_+(\gamma_1),\ldots,i_+(\gamma_{2g})$. 
(See \cite{fo} for this fact through 
the concept of {\it presentations of a pair of groups}.) 
By definition, $H_1 (M) / \Phi_+ =0$, and we have an exact sequence 
\[\begin{CD}
\Z^{2g+l} @>\sideset{^\mathfrak{t}}{}
{\Tmatrix{
\begin{smallmatrix} A \\ B \end{smallmatrix}}}\!\cdot >> 
\Z^{2g+l} @>>> H_1 (M) / \Phi_+ = 0.
\end{CD}\]
By the Hopfian property of $\Z^{2g+l}$, we see that 
$\sideset{^\mathfrak{t}}{}
{\Tmatrix{
\begin{smallmatrix} A \\ B \end{smallmatrix}}}$ is invertible. 
(1) follows from this. (See Remark \ref{acyclic}.)
 
(2) Through a standard argument using Eilenberg-MacLane spaces, 
we can assume that a given admissible presentation is obtained 
from a cell decomposition of $M$. Then 
$\left(\begin{smallmatrix} 
A \\ B \\C\end{smallmatrix}\right) \!\cdot$ gives the boundary map 
$C_2 (M,p;\KK_{N_k (M)}) \xrightarrow{\partial_2} 
C_1 (M,p;\KK_{N_k (M)})$. Considering 
the correspondence of 1-cycles, we have 
\[\begin{pmatrix} I_{2g} \\ 0_{(l,2g)} \\ 0_{2g} \end{pmatrix} 
-\begin{pmatrix} 0_{2g} \\ 0_{(l,2g)} \\ r'_k(M) \end{pmatrix} 
=\begin{pmatrix} A \\ B \\ C \end{pmatrix} X 
\in C_1 (M,p;\KK_{N_k (M)})\]
for some matrix $X$, where we write $0_{k}$ and $0_{(k,l)}$ 
for the zero matrices of sizes $k \times k$ and $k \times l$ 
respectively. (2) follows from this.
\end{proof}
\noindent
Note that from (2), we have $r'_k(M) = 
-C \begin{pmatrix} A \\ B \end{pmatrix}^{-1} \!
\begin{pmatrix} I_{2g} \\ 0_{(l,2g)}\end{pmatrix}$, namely the Magnus 
matrix $r'_k(M)$ can be computed from any 
admissible presentation of $\pi_1 (M)$. 

Next, we derive a formula for changing a generating system of 
$\pi_1 \Sg$. For a homology cylinder $(M,i_+,i_-)$, we take an 
admissible presentation of $\pi_1 M$ as in Definition \ref{admissible} 
and construct the matrices $A, B, C$ as before. 
Let $\gamma_1',\ldots,\gamma_{2g}'$ be 
another generating system of $\pi_1 \Sg$. We can take 
$\varphi \in \Aut \pi_1 \Sg$ such that 
$\gamma_i'=\varphi (\gamma_i)$ for $i=1,\ldots,2g$. 
\begin{prop}\label{basischange}
Let $r_k^\varphi(M)$ be the Magnus matrix corresponding to 
the new generating system. Then 
\[r_k^\varphi(M)=
\overline{\left(\frac{\partial \varphi(\gamma_j)}
{\partial \gamma_i}\right)}^{-1} r_k(M)
\sideset{^{\sigma_k(M)}\!}{}
{\mathop{\overline{
\left( \frac{\partial \varphi(\gamma_j)}{\partial \gamma_i} \right) 
}}\nolimits}.\]
\end{prop}
\begin{proof}
We have the following 
admissible presentation of $\pi_1 M$ with respect to 
$\gamma_1',\ldots,\gamma_{2g}'$:
\[\pi_1 M \cong 
\left\langle
\begin{array}{c|c}
\begin{array}{l}
i_- (\gamma_1'),\ldots,i_-(\gamma_{2g}'),\\
i_- (\gamma_1),\ldots,i_- (\gamma_{2g}),\\ 
z_1 ,\ldots, z_l, \\
i_+ (\gamma_1),\ldots,i_+ (\gamma_{2g}),\\
i_+ (\gamma_1'),\ldots,i_+ (\gamma_{2g}')
\end{array}&
\begin{array}{l}
i_- (\gamma_1')i_-(\varphi(\gamma_1))^{-1},\ldots, 
i_- (\gamma_{2g}')i_-(\varphi(\gamma_{2g}))^{-1}, \\
r_1, \ldots r_{2g+l},\\ 
i_+ (\gamma_1')i_+(\varphi(\gamma_1))^{-1},\ldots, 
i_+ (\gamma_{2g}')i_+(\varphi(\gamma_{2g}))^{-1},  
\end{array}
\end{array}\right\rangle.\]
The matrices $A',B',C'$ corresponding to this presentation 
are given by 
\begin{align*}
A' &= (I_{2g} \quad 0_{(2g,2g+l)} \quad 0_{2g})\\
B' &= 
\begin{pmatrix}
-\overline{\sideset{^{i_-}\!}{}{\Tmatrix{\frac{\partial 
\varphi(\gamma_j)}{\partial \gamma_i}}}}_{1 \le i,j \le 2g} & 
A & 0_{2g}\\
0_{(l,2g)} & B & 0_{(l,2g)}\\
0_{2g} & C & 
-\overline{\sideset{^{i_+}\!}{}{\Tmatrix{\frac{\partial 
\varphi(\gamma_j)}{\partial \gamma_i}}}}_{1 \le i,j \le 2g}
\end{pmatrix}\\
C' &= (0_{2g} \quad 0_{(2g,2g+l)} \quad I_{2g})
\end{align*}
\noindent
Computing $-C' \begin{pmatrix} A' \\ B' \end{pmatrix}^{-1}\!
\begin{pmatrix} I_{2g} \\ 0_{(4g+l,2g)}\end{pmatrix}$, 
we obtain the formula. 
\end{proof}

\vspace{10pt}

\section{Example: Relationship to the Gassner representation for 
string links}\label{sec:gassner}
In \cite{le}, Levine gave a method for constructing homology cylinders 
from pure string links. By this, we can obtain many 
homology cylinders not belonging to the subgroup $\Mg$. Also, we can 
see a relationship between the Gassner representation 
for string links and our representation. 

\begin{figure}[htbp]
\begin{center}
\includegraphics{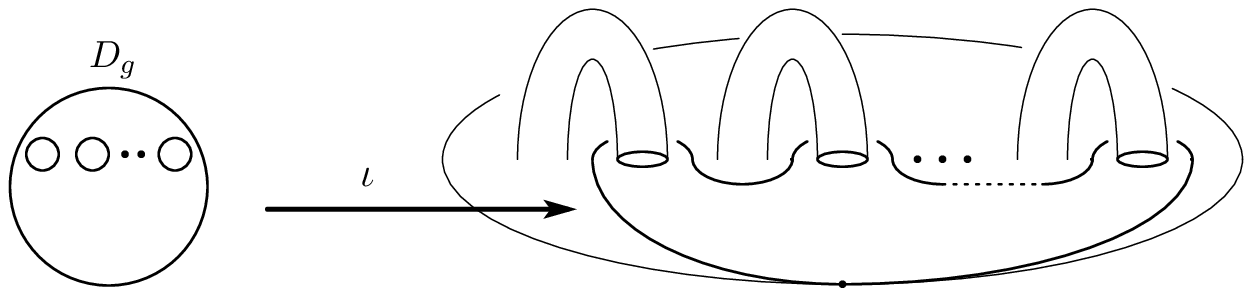}
\end{center}
\caption{}
\label{fig:ht}
\end{figure}

For a $g$-component pure string link $L \subset D^2 \times I$, 
we now construct a homology cylinder $M_L \in \Cg$ as follows. 
Consider a closed tubular neighborhood of 
the loops $\gamma_{g+1}, \gamma_{g+2}, \ldots, 
\gamma_{2g}$ in Figure 1 to be the image of 
an embedding $\iota:D_g \hookrightarrow \Sg$ of a $g$-holed disk $D_g$ 
as in Figure \ref{fig:ht}. 

\noindent
Let $C$ be the complement of an open tubular neighborhood of $L$ in 
$D^2 \times I$. For each choice of a framing of $L$, 
a homeomorphism $h:\partial C \xrightarrow{\cong} 
\partial (\iota(D_g) \times I)$ is fixed. 
Then the manifold $M_L$ obtained from $\Sg \times I$ 
by removing $\iota(D_g) \times I$ and regluing 
$C$ by $h$ becomes a homology cylinder. 
This construction gives an injective monoid homomorphism $\Lg \to \Cg$ 
from the monoid $\Lg$ of (framed) pure string links to $\Cg$. 
Moreover it also induces an injective homomorphism $\Slg \to \Hg$ from 
the concordance group $\Slg$ of (framed) pure string links to $\Hg$. 
In particular, the (smooth) knot concordance group, 
which coincides with $\mathcal{S}_{1}$, is embedded in $\Hg$. 
If we restrict these embeddings to the pure braid group, which is 
a subgroup of $\Lg$ and $\Slg$, their images are contained in $\Mg$. 

By the Gassner representation, we mean the crossed homomorphism 
$r_{G,k}:\Lg \rightarrow GL(g,$ $\KK_{N_k (D_g)})$ or 
$r_{G,k}:\Slg \rightarrow GL(g,\KK_{N_k (D_g)})$ 
given by a construction similar to that in the previous section. 
(In \cite{ld} and \cite{klw}, only $r_{G,2}$ is treated.) 
Comparing methods for calculating the Gassner and Magnus 
representations, we obtain the following. 
\begin{thm}\label{thm:conn}
For any pure string link $L \in \Lg$, 
$r_k (M_L) = 
\left(\begin{array}{cc}
\ast & 0_g \\ \ast & r_{G,k} (L)
\end{array}\right).$
\end{thm}
\noindent
We mention two remarks about this theorem. 
First we identify $F_g = \pi_1 D_g$ 
with the subgroup of $F_{2g} = \pi_1 \Sg$ generated by 
$\gamma_{g+1}, \ldots, \gamma_{2g}$. Then the maps 
$F_g = \langle \gamma_{g+1},\ldots,\gamma_{2g} \rangle$ 
$\hookrightarrow F_{2g} \twoheadrightarrow F_g$, where 
the second map sends $\gamma_1,\ldots, \gamma_g$ to $1$, show 
that $N_k (F_g) \subset N_k$ and 
$\KK_{N_k (F_g)} \subset \KK_{N_k}$. 
Second, the embeddings $\Lg \hookrightarrow \Cg$ and 
$\Slg \hookrightarrow \Hg$ have ambiguity with respect to 
framings. However we can check that 
the lower right part of $r_k (M_L)$ is independent of the framings. 
\begin{proof}[Proof of Theorem $\ref{thm:conn}$]
All we have to do is to give a suitable 
presentation of $\pi_1 M_L$. We divide $M_L$ into 
two parts $M$ and $C$ as follows. 

We take $g$ points $q_1,\ldots,q_g$ and 
$g$ paths $l_j$ from the base point $p$ to $q_j$ 
as in Figure \ref{fig:ht2}. 

\begin{figure}[htbp]
\begin{center}
\includegraphics{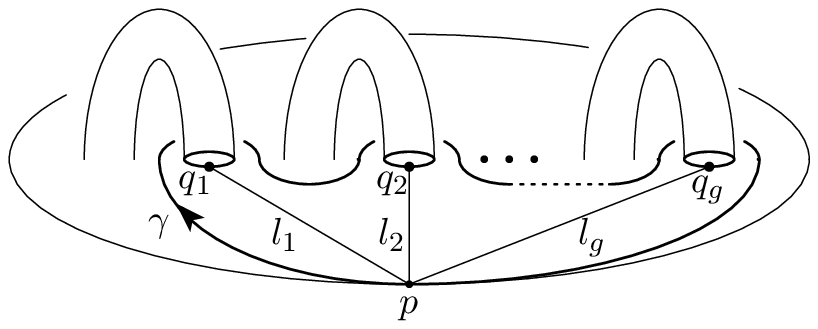}
\end{center}
\caption{}
\label{fig:ht2}
\end{figure}

\noindent
Let $M$ be the union of $\overline{\Sg \times I - 
D_g \times I}$ and $2g$ paths $i_+ (l_j)$ and $i_-(l_j)$ 
$(j=1,\ldots,g)$. Then 
\[\pi_1 M \cong 
\left\langle
\begin{array}{c|c}
\begin{array}{c}
i_- (\widehat{\gamma_1}),\ldots,
i_- (\widehat{\gamma_g})\\
i_+ (\widehat{\gamma_1}),\ldots,
i_+ (\widehat{\gamma_g})\\
i_+ (\gamma_{g+1}),\ldots,i_+ (\gamma_{2g}) \\
\delta_1,\ldots,\delta_g, i_+(\gamma) 
\end{array}
&
\begin{array}{c}
i_+ (\widehat{\gamma_1}) = i_-(\widehat{\gamma_1}) 
\delta_1 \\
\vdots \\
i_+ (\widehat{\gamma_g}) = i_-(\widehat{\gamma_g}) 
\delta_g
\end{array}
\end{array}\right\rangle\]
where $\widehat{\gamma_j}=[\gamma_1,\gamma_{g+1}] \cdots 
[\gamma_{j-1},\gamma_{g+j-1}] \gamma_j$, $\gamma$ is the loop 
corresponding to the outer boundary of $\iota(D_g)$ 
and $\delta_j$ is the composite of paths 
$i_- (l_j)$, $\overrightarrow{i_- (q_j) i_+ (q_j)}$ 
and $i_+ (l_j^{-1})$. 
We denote by $C$ the complement of 
an open tubular neighborhood of 
$L$ in $D_g \times I$ 
as before. 
\[\pi_1 C \cong \langle i_- (\gamma_{g+1}),\ldots, i_- (\gamma_{2g}), 
z_1,\ldots,z_l, i_+ (\gamma_{g+1}),\ldots, i_+ (\gamma_{2g}) \mid 
r_1, \ldots r_{g+l}
\rangle\]
is given by the Wirtinger presentation of 
$D \times I - L$. 
We glue $C$ to $M$ by using some fixed framing. 
Then it is easily seen that 
$\pi_1 (M \cap C)$ is the free group generated by 
$\{i_+ (\gamma_{g+1}),\ldots, i_+ (\gamma_{2g}),$ 
$\delta_1,\ldots,\delta_g, i_+(\gamma)\}$.

Using the above decomposition, we obtain
\[\pi_1 M_L \cong 
\left\langle
\begin{array}{c|c}
\begin{array}{c}
i_- (\widehat{\gamma_1}),\ldots,
i_- (\widehat{\gamma_g})\\
i_- (\gamma_{g+1}),\ldots,i_- (\gamma_{2g}) \\
z_1 ,\ldots, z_{l}\\ 
i_+ (\widehat{\gamma_1}),\ldots,
i_+ (\widehat{\gamma_g})\\
i_+ (\gamma_{g+1}),\ldots,i_+ (\gamma_{2g}) 
\end{array}
&
\begin{array}{c}
i_+ (\widehat{\gamma_1}) = i_-(\widehat{\gamma_1}) 
\widehat{\delta_1} \\
\vdots \\
i_+ (\widehat{\gamma_g}) = i_-(\widehat{\gamma_g}) 
\widehat{\delta_g}
\\
r_1, \ldots r_{g+l}
\end{array}
\end{array}\right\rangle\]
where $\widehat{\delta_i}$ are words 
in $i_- (\gamma_{g+1}),\ldots, i_- (\gamma_{2g}), 
z_1,\ldots,z_l$, 
$i_+ (\gamma_{g+1}),\ldots, i_+ (\gamma_{2g})$ which 
depend on the framing. 
Rewrite the 
above presentation by using $i_+ (\gamma_j)$'s 
and $i_- (\gamma_j)$'s instead of $i_+ (\widehat{\gamma_j})$'s 
and $i_- (\widehat{\gamma_j})$'s. 
This process does not affect generators 
$i_- (\gamma_{g+j}), z_j, i_+ (\gamma_{g+j})$ and relations $r_j$. 
From the resulting admissible presentation, 
we can compute the Magnus matrix of $M_L$. 
Then our claim follows by comparing it with the method for 
calculating the Gassner matrix of $L$ from the Wirtinger presentation 
of $\pi_1 C$, which is given in \cite{klw}.
\end{proof}
\begin{cor}\label{notnormal}
$\Mg$ is not a normal subgroup of $\Hg$ for $g \ge 3$.
\end{cor}
\begin{proof}
In \cite{klw}, they gave 3-component pure string links denoted by 
$L_5$ and $L_6$ having the condition 
that $L_5$ is a pure braid, while the 
conjugate $L_6 L_5 L_6^{-1}$ is not. To show that 
$L_6 L_5 L_6^{-1}$ is not a pure braid, they use the 
fact that $r_{G,2}(L_6 L_5 L_6^{-1})$ has an entry not 
belonging to $\Z N_2 (D_3)$. 
Then our claim follows 
from Theorem \ref{thm:conn} with respect to this example.
\end{proof}
\begin{example}\label{eg4}
Let $L$ be a 2-component pure string link as depicted in 
Figure \ref{fig:string}. 

\begin{figure}[htbp]
\begin{center}
\includegraphics{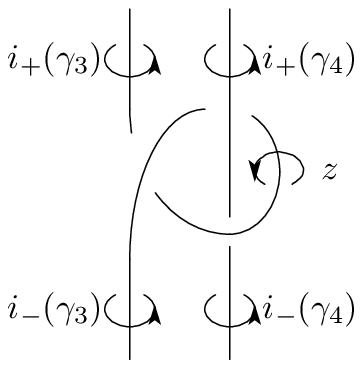}
\end{center}
\caption{}
\label{fig:string}
\end{figure}

\noindent
Then the presentation of $\pi_1 M_L$ 
is given by 
{\small 
\[\pi_1 M_L \cong 
\left\langle
\begin{array}{c|c}
\begin{array}{c}
i_-(\gamma_1),\ldots,i_-(\gamma_4)\\
z\\
i_+(\gamma_1),\ldots,i_+(\gamma_4)\\
\end{array} &
\begin{array}{l}
i_+(\gamma_1) i_-(\gamma_3)^{-1} i_+(\gamma_4) i_-(\gamma_1)^{-1},\\
{[i_+(\gamma_1),i_+(\gamma_3)]} i_+(\gamma_2) z i_-(\gamma_2)^{-1} 
{[i_-(\gamma_3),i_-(\gamma_1)]},\\
i_+(\gamma_4) i_-(\gamma_3) i_+(\gamma_4)^{-1} z^{-1},\quad
i_-(\gamma_3) i_+(\gamma_3)^{-1} i_-(\gamma_3)^{-1} z,\\
i_-(\gamma_4) z^{-1} i_+(\gamma_4)^{-1} z,\quad
\end{array}
\end{array}\right\rangle,\]
}
\noindent
where we use the blackboard framing. 
We now compute $r_2(M_L)$. 
We identify $N_2$ and $N_2 (M_L)$ by 
using $i_+$. Using the presentation, we have $z=i_-(\gamma_3)=\gamma_3$, 
$i_-(\gamma_4)=\gamma_4$, $i_-(\gamma_2)=\gamma_2\gamma_3$ and 
$i_-(\gamma_1)=\gamma_1\gamma_3^{-1}\gamma_4$ in $N_2$. Then 
\begin{align*}
\begin{pmatrix}A \\ B \end{pmatrix} &=
\begin{pmatrix}
-1&\gamma_3^{-1}-1& 0 & 0 & 0 \\
0 & -1 & 0 & 0 & 0 \\
-\gamma_1^{-1}\gamma_3 & 
1- \gamma_1^{-1} \gamma_3 \gamma_4^{-1} &
\gamma_4^{-1}& 1-\gamma_3 & 0\\
0 & 0 & 0 & 0 & 1 \\
0 & \gamma_2^{-1}& -1 & \gamma_3 & \gamma_3-\gamma_3 \gamma_4^{-1}
\end{pmatrix}, \\
C &=\begin{pmatrix}
1& 1-\gamma_3^{-1} &0&0&0\\
0&1&0&0&0\\
0&\gamma_1^{-1} -1&0&-1&0\\
\gamma_1^{-1}\gamma_3&0&1-\gamma_3^{-1}& 0&-\gamma_3 
\end{pmatrix}.
\end{align*}
Hence
\[r_2(M_L) =
\left(\begin{array}{cccc}
1&0&0&0\\
&&&\\
0&1&0&0
\\
&&&\\
\frac{-\gamma_1^{-1}}{\gamma_3^{-1}+\gamma_4^{-1}-1}
&\frac{\gamma_2^{-1} \gamma_3^{-1} \gamma_4^{-1} 
-\gamma_4^{-1}+1}{\gamma_3^{-1}+\gamma_4^{-1}-1} &
\frac{\gamma_3^{-1}}{\gamma_3^{-1}+\gamma_4^{-1}-1}
&\frac{\gamma_4^{-1}(\gamma_4^{-1}-1)}{\gamma_3^{-1}+\gamma_4^{-1}-1}\\
&&&\\
\frac{\gamma_1^{-1}\gamma_3\gamma_4^{-1}}{\gamma_3^{-1}+\gamma_4^{-1}-1}
& \frac{(1-\gamma_3^{-1})(\gamma_2^{-1} \gamma_3^{-1} 
-\gamma_2^{-1}-1)}{\gamma_3^{-1}+\gamma_4^{-1}-1} 
&\frac{\gamma_3^{-1}-1}{\gamma_3^{-1}+\gamma_4^{-1}-1}
&\frac{-\gamma_3^{-1}\gamma_4^{-1}+\gamma_3^{-1}
+2\gamma_4^{-1}-1}{\gamma_3^{-1}+\gamma_4^{-1}-1}
\end{array}\right).\]
\end{example}

\vspace{10pt}

\section{Higher-order Alexander invariants and 
torsion-degree functions}\label{sec:harvey}
Here we summarize the theory of 
higher-order Alexander invariants along the lines of Harvey's 
papers \cite{har,har2}. 
For our use, we generalize them to functions of matrices 
called {\it torsion-degree functions}. 

A group $\Gamma$ is {\it poly-torsion-free-abelian} (PTFA, for short) 
if $\Gamma$ has a normal series of finite length whose 
successive quotients are all torsion-free abelian. In particular, 
free nilpotent quotients $N_k$ are PTFA for all $k \ge 2$. 
Note that every subgroup of 
a PTFA group is also PTFA. 
For each PTFA group $\Gamma$, the group ring $\Z \Gamma$ is known to 
be an Ore domain, 
so that it can be embedded in the {\it right field of fractions}
$\KK_{\Gamma}:=\Z \Gamma (\Z \Gamma -\{ 0 \})^{-1}$, which is a skew field. 
We refer to \cite{co2}, \cite{pa} for localizations 
of non-commutative rings. 

We will also use the following localizations of $\Z \Gamma$ placed 
between $\Z \Gamma$ and $\KK_\Gamma$. 
Let $\psi : \Gamma \twoheadrightarrow \Z$ be an epimorphism. 
Then we have an exact sequence 
\[1 \longrightarrow (\Gamma^\psi:= \Ker \psi) \longrightarrow 
\Gamma \stackrel{\psi}{\longrightarrow} \Z \longrightarrow 1.\]
We take a splitting $\xi:\Z \to \Gamma$ and 
put $\widetilde{t}:=\xi(1) \in \Gamma$. 
Since $\Gamma^\psi$ is again 
a PTFA group, 
$\Z \Gamma^\psi$ can be embedded in 
its right field of fractions 
$\KK_{\Gamma^\psi}=\Z \Gamma^\psi (\Z \Gamma^\psi -\{ 0 \})^{-1}$. 
Moreover, we can also construct 
$\Z \Gamma (\Z \Gamma^\psi -\{ 0 \})^{-1}$. 
Then the splitting $\xi$ gives 
an isomorphism between $\Z \Gamma (\Z \Gamma^\psi -\{ 0 \})^{-1}$ and 
the skew Laurent polynomial ring $\KK_{\Gamma^\psi} [t^{\pm}]$, 
in which $at=t (\widetilde{t}^{-1} a \widetilde{t}\,)$ holds 
for each $a \in \KK_{\Gamma^\psi}$. 
$\KK_{\Gamma^\psi} [t^{\pm}]$ is known to be a non-commutative 
right and left principal ideal domain. 
By definition, we have inclusions 
\[\Z \Gamma \hookrightarrow \KK_{\Gamma^\psi} [t^{\pm}] 
\hookrightarrow \KK_{\Gamma}.\]
$\KK_{\Gamma^\psi} [t^{\pm}]$ and $\KK_{\Gamma}$ are known to be flat 
$\Z \Gamma$-modules. 
On $\KK_{\Gamma^\psi} [t^{\pm}]$, we have a map 
$\deg^{\psi}:\KK_{\Gamma^\psi} [t^{\pm}] \to \Z_{\ge 0} \cup \{ \infty \}$ 
assigning to each polynomial its degree. 
We put $\deg^{\psi}(0):= \infty $. By this, 
$\KK_{\Gamma^\psi} [t^{\pm}]$ becomes a Euclidean domain. 
The composite $\Z \Gamma (\Z \Gamma^\psi -\{ 0 \})^{-1} 
\xrightarrow{\cong} 
\KK_{\Gamma^{\psi}} [t^{\pm}] \xrightarrow{\deg^{\psi}} 
\Z_{\ge 0} \cup \{ \infty \}$ does not depend on 
the choice of $\xi$. 

Harvey's higher-order Alexander invariants \cite{har2} 
are defined as follows. 
Let $G$ be a finitely presentable group 
and let $\varphi:G \twoheadrightarrow \Z$ be 
an epimorphism. For a PTFA group $\Gamma$ and an epimorphism 
$\varphi_{\Gamma} : G \twoheadrightarrow \Gamma$, 
$(\varphi_{\Gamma},\varphi)$ is 
called an {\it admissible pair} for $G$ if there exists an epimorphism 
$\psi:\Gamma \twoheadrightarrow \Z$ satisfying 
$\varphi=\psi \circ \varphi_{\Gamma}$. 
For each admissible pair $(\varphi_{\Gamma},\varphi)$ for $G$, 
we regard $\KK_{\Gamma^\psi}[t^{\pm}]=
\Z \Gamma (\Z \Gamma^\psi -\{ 0 \})^{-1}$ as a $\Z G$-module, and 
we define the higher-order Alexander invariant for 
$(\varphi_{\Gamma},\varphi)$ by 
\begin{align*}
\overline{\delta}_\Gamma^\psi (G) &= \dim_{\KK_{\Gamma^\psi}} 
\bigl( H_1 (G;\KK_{\Gamma^\psi}[t^{\pm}]) \bigr) \in \Z_{\ge 0} 
\cup \{ \infty \},\\
\delta_\Gamma^\psi (G) &= \dim_{\KK_{\Gamma^\psi}} 
\bigl( T_{\KK_{\Gamma^\psi}[t^{\pm}]} H_1 (G;\KK_{\Gamma^\psi}[t^{\pm}]) 
\bigr) \in \Z_{\ge 0},
\end{align*}
where $T_{\KK_{\Gamma^\psi}[t^{\pm}]} M$ denotes the 
$\KK_{\Gamma^\psi}[t^{\pm}]$-torsion part of a 
$\KK_{\Gamma^\psi}[t^{\pm}]$-module $M$. 
$\overline{\delta}_\Gamma^\psi (G)$ and $\delta_\Gamma^\psi (G)$ are 
called the $\Gamma$-{\it degree}
and the {\it refined $\Gamma$-degree} respectively. 
(Our definition is slightly different from 
that of \cite{har2}.) 
Note that the right $\KK_{\Gamma^\psi}[t^{\pm}]$-module 
$H_1 (G;\KK_{\Gamma^\psi}[t^{\pm}])$ can be decomposed into 
\[H_1 (G;\KK_{\Gamma^\psi}[t^{\pm}]) = 
(\KK_{\Gamma^\psi}[t^{\pm}])^r \oplus \left(\bigoplus_{i=1}^l 
\frac{\KK_{\Gamma^\psi}[t^{\pm}]}{p_i (t) 
\KK_{\Gamma^\psi}[t^{\pm}]}\right)\]
for some $r \in \Z_{\ge 0}$ and $p_i (t) \in \KK_{\Gamma^\psi}[t^{\pm}]$, 
then 
\begin{align*}
\overline{\delta}_\Gamma^\psi (G) &= 
\begin{cases}
\sum_{i=1}^l \deg^\psi (p_i (t)) & (r=0),\\
\infty & (r > 0)
\end{cases},\\
\delta_\Gamma^\psi (G) &= 
\sum_{i=1}^l \deg^\psi (p_i (t)).
\end{align*}
For a connected space $X$ and an admissible pair 
$(\varphi_{\Gamma},\varphi)$ for $\pi_1 X$, we define 
$\overline{\delta}_\Gamma^\psi (X):= 
\overline{\delta}_\Gamma^\psi (\pi_1 X)$ and 
$\delta_\Gamma^\psi (X):=\delta_\Gamma^\psi (\pi_1 X)$.

For a finitely presentable group $G$ and an admissible pair 
$(\varphi_{\Gamma},\varphi)$ for $G$, the (refined) 
$\Gamma$-degree can be computed from 
any presentation matrix of 
the right $\KK_{\Gamma^\psi}[t^{\pm}]$-module 
$H_1 (G;\KK_{\Gamma^\psi}[t^{\pm}])$. Therefore we can 
consider it to be a function on the set 
$M(\KK_{\Gamma^\psi}[t^{\pm}])$ of all 
matrices with entries in $\KK_{\Gamma^\psi}[t^{\pm}]$. 
Here we extend this function to $M (\KK_\Gamma)$ as follows. 

First, we extend $\deg^\psi:\KK_{\Gamma^\psi}[t^{\pm}] \to 
\Z_{\ge 0} \cup \{ \infty \}$ to $\deg^\psi:\KK_{\Gamma} \to 
\Z \cup \{ \infty \}$ by setting 
$\deg^\psi (f g^{-1})= \deg^\psi(f) - \deg^\psi(g)$ 
for $f \in \Z \Gamma, g \in \Z \Gamma - \{ 0 \}$ 
(see for instance \cite[Proposition 9.1.1]{co2}). 
It induces a group homomorphism 
$\deg^{\psi}:(\KK_{\Gamma}^\times)_{\ab} \to \Z$, 
where $(\KK_{\Gamma}^\times)_{\ab}$ 
is the abelianization of the multiplicative group 
$\KK_{\Gamma}^\times = \KK_{\Gamma} - \{ 0 \}$. 

Since $\KK_{\Gamma}$ is a skew field, 
we have the Dieudonn\'e determinant 
\[\det : GL(\KK_\Gamma) \longrightarrow 
(\KK_\Gamma^{\times})_{\mathrm{ab}},\]
which is a homomorphism characterized by 
the following three properties:
\begin{itemize}
\item[(a)] $\det I=1$.
\item[(b)] If $A'$ is obtained by multiplying 
a row of a matrix $A \in GL(\KK_\Gamma)$ 
by $a \in \KK_\Gamma^\times$ from the left, 
then $\det A'= a \cdot \det A$. 
\item[(c)] If $A'$ is obtained by adding to a row of 
a matrix $A$ 
a left $\KK_\Gamma$-linear combination of other rows, 
then $\det A'= \det A$.
\end{itemize}
It induces 
an isomorphism between $K_1 (\KK_\Gamma) \xrightarrow{\cong} 
(\KK_\Gamma^\times)_{\mathrm{ab}}$. 

The following lemma will be used in our generalization 
of Harvey's invariants. We denote by $M(m,n,\KK_{\Gamma})$ 
the set of all $m \times n$ matrices with entries in $\KK_{\Gamma}$. 
\begin{lem}\label{lem:torsion}
For $A \in M(m,n,\KK_{\Gamma})$ with $\rank_{\KK_{\Gamma}} A=k$, 
let $U \in M(m-k,m,\KK_{\Gamma})$, $V \in M(n,n-k,\KK_{\Gamma})$ be 
matrices satisfying 
\[\begin{cases}
UA = 0, & \rank_{\KK_{\Gamma}} U=m-k, \\
AV = 0, & \rank_{\KK_{\Gamma}} V=n-k.
\end{cases}\]
For each $I \subset \{ 1,2,\ldots,m\}$, $J \subset \{ 1,2,\ldots,n\}$ 
with $\# I=m-k$, $\# J=n-k$, let $U_I$ denote the square matrix defined 
by taking $i$-th columns from $U$ for all $i \in I$, and $V_J$ denote 
the one defined by taking $j$-th rows from $V$ 
for all $j \in J$. We also denote by $A_{I^c J^c}$ 
the one defined 
by taking $i$-th rows from $A$ for all 
$i \in I^c:=\{ 1,2,\ldots,m\} - I$ and then taking 
$j$-th columns for all $j \in J^c:=\{ 1,2,\ldots,n\} - J$. 

\noindent
$(1)$ If $U_I$ or $V_J$ is not invertible, 
then $A_{I^c J^c}$ is not invertible.

\noindent
$(2)$ Otherwise, 
\[\Delta(A;U,V):=\sgn(I I^c) \sgn(J J^c) 
\frac{\det A_{I^c J^c}}{\det U_I \det V_J} 
\in (\KK_{\Gamma}^{\times})_{\ab}\] 
is independent of the choice of $I$ and $J$ such that 
$U_I$, $V_J$ are invertible, where $\sgn(I I^c) \in \{ \pm 1 \}$ 
$($resp.~$\sgn(J J^c))$ 
is the signature of the juxtaposition of $I$ and $I^c$ 
$($resp.~$J$ and $J^c)$, 
and we put $\det \emptyset :=1$.

\noindent
$(3)$ For $P_1 \in GL(m,\KK_{\Gamma})$, $P_2 \in GL(n,\KK_{\Gamma})$, 
$Q_1 \in GL(m-k,\KK_{\Gamma})$ and $Q_2 \in GL(n-k,\KK_{\Gamma})$, 
\[\Delta(P_1^{-1}AP_2^{-1};Q_1 U P_1,P_2 V Q_2)=
\frac{\Delta(A;U,V)}{\det P_1 \det P_2 \det Q_1 \det Q_2}.\]
\end{lem}
\begin{proof}
(1) and (2) are deduced from easy observation using non-commutative 
linear algebra. 
To prove (3), it suffices to show in the cases where 
$P_1,P_2,Q_1,Q_2$ are matrices of elementary transformations, 
and it can be easily checked. 
\end{proof}
\begin{remark}
In the above situation, the sequence 
\[\begin{CD}
0 @>>> \KK_{\Gamma}^{n-k} @>{V \cdot}>> \KK_{\Gamma}^n @>{A \cdot}>>
\KK_{\Gamma}^m @>{U \cdot}>>\KK_{\Gamma}^{m-k} @>>> 0
\end{CD}\]
is exact. By taking the standard basis for each vector space, 
we regard the sequence as a based acyclic chain complex. 
Then we can take its torsion 
(see \cite{mi}, \cite{tu2} for generalities of torsions). 
This torsion coincides with $\Delta(A;U,V)$ up to sign.
\end{remark}

As seen in Lemma \ref{lem:torsion} (3), 
$\Delta(A;U,V)$ does depend on 
$U$ and $V$. The following definition and lemma give 
our rule to take $U$ and $V$. 
\begin{definition}\label{def:primitive}
Let $A \in M(m,n,\KK_{\Gamma})$ with $\rank_{\KK_\Gamma} A =k$. 
$(U,V)$ is said to be $\psi$-{\it primitive} for $A$ if 
\begin{itemize}
\item[(1)] $U$, $V$ have 
entries in $\KK_{\Gamma^\psi} [t^\pm]$. 
\item[(2)] The row vectors $u_1,\ldots,u_{m-k} \in 
(\KK_{\Gamma^\psi} [t^\pm])_m$ of $U$ generate 
$\Ker (\cdot A) \cap (\KK_{\Gamma^\psi} [t^\pm])_m$ in $(\KK_\Gamma)_m$ 
as a left $\KK_{\Gamma^\psi} [t^\pm]$-module.
\item[(3)] The column vectors $v_1,\ldots,v_{n-k} \in 
(\KK_{\Gamma^\psi} [t^\pm])^n$ of $V$ generate 
$\Ker (A \cdot) \cap (\KK_{\Gamma^\psi} [t^\pm])^n$ in $(\KK_\Gamma)^n$ 
as a right $\KK_{\Gamma^\psi} [t^\pm]$-module.
\end{itemize}
\end{definition}
\begin{lem}\label{lem:UV}
$(1)$ There exists a pair $(U,V)$ which is $\psi$-primitive for $A$. \\
$(2)$ Suppose $U \in M(m-k,m,\KK_{\Gamma^\psi} [t^\pm])$ and 
$V \in M(n,n-k,\KK_{\Gamma^\psi} [t^\pm])$ satisfy $UA=0$ and $AV=0$. 
$(U,V)$ is $\psi$-primitive for $A$ 
if and only if there exist $\widetilde{P}_1 
\in GL(m,\KK_{\Gamma^\psi} [t^\pm])$ and 
$\widetilde{P}_2 \in GL(n,\KK_{\Gamma^\psi} [t^\pm])$ satisfying 
\[U \widetilde{P}_1 = (0_{(m-k,k)} \quad I_{m-k}), \quad 
\widetilde{P}_2 V = (0_{(n-k,k)} \quad I_{n-k})^T.\]
$(3)$ If $(U,V)$ and $(U',V')$ are $\psi$-primitive for $A$, 
then there exist 
$P_1 \in GL(m,\KK_{\Gamma^\psi} [t^\pm])$, 
$P_2 \in GL(n,$ $\KK_{\Gamma^\psi} [t^\pm])$, 
$Q_1 \in GL(m-k,\KK_{\Gamma^\psi} [t^\pm])$ and 
$Q_2 \in GL(n-k,\KK_{\Gamma^\psi} [t^\pm])$ such that
\[U P_1 = U',\quad P_2 V = V',\quad Q_1 U = U',\quad V Q_2 =V'.\]
\end{lem}
\begin{proof}
We prove only for $V$. \\
(1) For right 
$\KK_{\Gamma^\psi} [t^\pm]$-homomorphisms $(\KK_{\Gamma^\psi} [t^\pm])^n 
\stackrel{i}{\hookrightarrow} \KK_\Gamma^n \xrightarrow{A\cdot} 
\KK_\Gamma^m$, $\Ker ((A\cdot)\circ i) 
= \Ker (A\cdot) \cap (\KK_{\Gamma^\psi} [t^\pm])^n$ is a right 
free $\KK_{\Gamma^\psi} [t^\pm]$-module of rank $n-k$. We take a 
basis $v_1,\ldots,v_{n-k}$ and put $V=(v_1,\ldots,v_{n-k})$. 
Then $V$ satisfies the desired property. \\
(2) Suppose $V$ generates $\Ker (A\cdot) \cap 
(\KK_{\Gamma^\psi} [t^\pm])^n$. The quotient module 
$(\KK_{\Gamma^\psi} [t^\pm])^n/\Ker ((A\cdot)\circ i)$ is 
$\KK_{\Gamma^\psi} [t^\pm]$-torsion free, and hence 
$\KK_{\Gamma^\psi} [t^\pm]$-free. Taking a splitting, 
we have a direct sum decomposition $(\KK_{\Gamma^\psi} [t^\pm])^n \cong 
\bigl( (\KK_{\Gamma^\psi} [t^\pm])^n/\Ker ((A\cdot)\circ i)\bigr) 
\oplus \Ker ((A\cdot)\circ i)$. 
We can extend $V$ to obtain a basis 
$(\tilde{v_1}, \ldots, \tilde{v_k}, V)$ for 
$(\KK_{\Gamma^\psi} [t^\pm])^n$. Then 
$\widetilde{P}_2:=(\tilde{v_1}, \ldots, \tilde{v_k}, V)^{-1}$ 
satisfies $\widetilde{P}_2 V = (0_{(n-k,k)} \quad I_{n-k})^T$. 
The inverse is clear. \\
(3) The existence of $P_2$ follows immediately from (2). 
That of $Q_2$ is also clear since $V$ and $V'$ are 
bases of the same right $\KK_{\Gamma^\psi} [t^\pm]$-module.
\end{proof}

\begin{definition}\label{def:TFfunc}
Let $\Gamma$ be a PTFA group 
and let $\psi:\Gamma \twoheadrightarrow \Z$ be an 
epimorphism. \\
(1) The {\it torsion-degree function} 
$d_\Gamma^\psi : M (\KK_\Gamma) \to \Z$ is defined by 
\[d_\Gamma^\psi (A) := \deg^\psi (\Delta(A;U,V))\]
for a pair $(U,V)$ which is $\psi$-primitive for $A$.\\
(2) The {\it truncated torsion-degree function} 
$\overline{d}_\Gamma^\psi : M (\KK_\Gamma) \to \Z \cup \{\infty\}$ 
is defined by 
\[\overline{d}_\Gamma^\psi (A) := 
\begin{cases}
d_\Gamma^\psi (A) & \text{if $\rank A \ge m-1$,}\\
\infty & \text{otherwise}
\end{cases}\]
for $A \in M(m,n,\KK_\Gamma)$.
\end{definition}
\noindent
Since $\KK_{\Gamma^\psi} [t^\pm]$ is a Euclidean domain, 
every $P \in GL(\KK_{\Gamma^\psi} [t^\pm])$ can be decomposed as 
products of elementary matrices and diagonal matrices in 
$GL(\KK_{\Gamma^\psi} [t^\pm])$, which 
shows that $\deg^\psi (\det P)=0$. 
Lemmas \ref{lem:torsion} and \ref{lem:UV} together with this fact 
show that $d_\Gamma^\psi$ and 
$\overline{d}_\Gamma^\psi$ are well-defined. 
\begin{example}\label{ex:TDF}
(1) For $A \in GL(\KK_\Gamma)$, we have 
$d_\Gamma^\psi (A) =\overline{d}_\Gamma^\psi (A)= \deg^\psi (\det A)$.\\
(2) Let $M$ be a finitely generated right 
$\KK_{\Gamma^\psi} [t^\pm]$-module, and let $A$ be a presentation 
matrix of $M$. Then we have 
$d_\Gamma^\psi (A) = \dim_{\KK_{\Gamma^\psi}}
\bigl(T_{\KK_{\Gamma^\psi} [t^\pm]} M \bigr)$. As for 
$\overline{d}_\Gamma^\psi (A)$, we can see that 
$\overline{d}_\Gamma^\psi (A) \in \Z$ 
if and only if 
the rank of the $\KK_{\Gamma^\psi} [t^\pm]$-free part of $M$ is 
less than $2$.\\
(3) Let $G$ be a finitely presentable group. We take a presentation 
$\langle x_1,\ldots,x_l \mid r_1,\ldots,r_m \rangle$ of $G$. For 
an admissible pair $(\varphi_\Gamma,\varphi)$, 
the matrix $A:=
\overline{\sideset{^{\varphi_\Gamma}\!}{}
{\Tmatrix{\frac{\partial r_j}{\partial x_i}}}}_{\begin{subarray}{c}
{}1 \le i \le l\\
1 \le j \le m
\end{subarray}}$ at 
$\KK_{\Gamma^\psi} [t^\pm]$ gives a presentation matrix of 
$H_1 (G,\{ 1 \};\KK_{\Gamma^\psi} [t^\pm])$. 
Then Harvey's invariants are given by 
\begin{align*}
\delta_\Gamma^\psi (G) &= \dim_{\KK_{\Gamma^\psi}} 
\bigl( T_{\KK_{\Gamma^\psi}[t^{\pm}]} H_1 (G;\KK_{\Gamma^\psi}[t^{\pm}]) 
\bigr) = d_\Gamma^\psi (A),\\
\overline{\delta}_\Gamma^\psi (G) &= \dim_{\KK_{\Gamma^\psi}} 
\bigl( H_1 (G;\KK_{\Gamma^\psi}[t^{\pm}]) \bigr) = 
\overline{d}_\Gamma^\psi (A),
\end{align*}
where the second equality of each case follows from 
the direct sum decomposition 
\[H_1 (G,\{1\};\KK_{\Gamma^\psi}[t^{\pm}]) \cong 
H_1 (G;\KK_{\Gamma^\psi}[t^{\pm}]) \oplus \KK_{\Gamma^\psi}[t^{\pm}]\]
shown by Harvey in \cite{har}.
\end{example}
\begin{remark}
Friedl \cite{fr} gave an interpretation of Harvey's invariants by 
Reidemeister torsions. The definition of our truncated torsion-degree 
functions has some overlaps with his description. 
\end{remark}

\vspace{10pt}

\section{Applications of torsion-degree functions to 
homology cylinders}\label{sec:apply}
In this section, we study some invariants of homology cylinders 
arising from the Magnus representation and Reidemeister torsions 
by using torsion-degree functions associated 
to nilpotent quotients $N_k$ of $\pi_1 \Sg$. 
$N_k$ is known to be a finitely generated torsion-free nilpotent group. 
In particular, it is PTFA. 

Note that we can take a primitive element of $H^1 (\Sg)$ instead of 
an epimorphism $N_k \twoheadrightarrow \Z$ to define 
a torsion-degree function 
since $\Hom (N_k,\Z)=H^1 (N_k) = H^1 (N_2) = H^1 (\Sg)$. 
We denote by $PH_1 (\Sg)$ the set of primitive elements of $H^1 (\Sg)$. 

\vspace{8pt}

\subsection{The Magnus representation and torsion-degree functions}
\label{subsec:Magnus}
First, we apply torsion-degree functions to the Magnus matrix. However, 
it turns out that the result is trivial. 
\begin{thm}\label{vanish}
Let $M$ be a homology cylinder. For any $\psi \in PH^1 (\Sg)$, 
the torsion-degree $d_{N_k}^\psi (r_k (M))$ is always zero. 
\end{thm}
\begin{proof}
By definition, $d_{N_k}^\psi$ is additive for products of invertible 
matrices, and invariant under taking the transpose and 
operating the involution. 
Moreover, it vanishes for matrices in $GL(\Z N_k)$. 
In \cite{sa5}, we show that there exists a matrix 
$\widetilde{J} \in GL(2g,\Z N_k)$ satisfying the equality 
\[\overline{r_k (M)^T}\ \widetilde{J}\ r_k (M) 
={}^{\sigma_k(M)}\widetilde{J}.\]
By applying $d_{N_k}^\psi$ to it, we obtain 
$2d_{N_k}^\psi (r_k (M))=0$. 
This completes the proof.
\end{proof}
\begin{example}
Consider the homology cylinder $M_L$ in Example \ref{eg4}. 
$d_{N_2}^\psi (r_2 (M_L))$ is given by the 
degree of $\det r_2 (M_L) =\frac{\gamma_3+\gamma_4-1}
{\gamma_3 \gamma_4(\gamma_3^{-1}+\gamma_4^{-1}-1)} $ 
with respect to $\psi$. It is zero. 
\end{example}

To extract some numerical information from $r_k(M)$, we next 
apply torsion-degree functions to $I_{2g}-r_k(M)$. 
Here we assume $M \in \Cg [k]$ 
and consider only $\overline{d}_{N_k}^\psi$. 
The function $\overline{d}_{N_k}^\psi 
(I_{2g}-r_k(\cdot)):\Cg [k] \to \Z \cup \{\infty\}$ 
factors through $\Hg$ since $r_k$ does. 
Note that for every $(M,i_+,i_-) \in \Cg [k]$, 
two inclusions $i_+$ and $i_-$ induce the same isomorphism 
$i_+=i_-:N_k \xrightarrow{\cong} N_k (M)$,
so that we can naturally identify them. 
Under this identification, we have the following. 

\begin{lem}\label{uandw}
Let $M$ be a homology cylinder belonging to $\Cg [k]$.
\begin{itemize}
\item[$(1)$] $(1-\gamma_1^{-1}, \ldots ,1-\gamma_{2g}^{-1})
(I_{2g}-r_k(M))=0$.
\item[$(2)$] $(I_{2g}-r_k(M))\overline{\left( 
\frac{\partial \zeta}{\partial \gamma_1},\ldots,
\frac{\partial \zeta}{\partial \gamma_{2g}}\right)}^T
=0$.
\end{itemize}
\end{lem}
\begin{proof}
We take an admissible presentation of $\pi_1 M$ 
as in Definition \ref{admissible}. We also take 
the matrices $A, B, C \in \Z N_k$ corresponding to it. 
For simplicity, we put $\overrightarrow{1-\overline{\gamma}}
:=(1-\gamma_1^{-1},\ldots,1-\gamma_{2g}^{-1} )$, 
$\overrightarrow{1-\overline{z}}
:=(1-z_1^{-1},\ldots,1-z_{l}^{-1})$ 
and 
$\zegamma
:=\left( 
\frac{\partial \zeta}{\partial \gamma_1},\ldots,
\frac{\partial \zeta}{\partial \gamma_{2g}}\right)$.

(1) Using Fundamental formula of free calculus 
(see \cite[Proposition 3.4]{bi}), we have 
\[(\overrightarrow{1-\overline{\gamma}} \quad 
\overrightarrow{1-\overline{z}} 
\quad \overrightarrow{1-\overline{\gamma}})
\begin{pmatrix} A \\ B \\C\end{pmatrix}=0.\]
Then, by Proposition \ref{howto}, 
\[(\overrightarrow{1-\overline{\gamma}} \quad 
\overrightarrow{1-\overline{z}})=-(\overrightarrow{1-\overline{\gamma}})
C\begin{pmatrix} A \\ B \end{pmatrix}^{-1}
=(\overrightarrow{1-\overline{\gamma}})(r_k(M) \quad Z).\]
Our claim follows by taking their first 2g columns. 

(2) Let $\tau_{\zeta} \in \Mg \subset \Cg$ be 
the Dehn twist along $\zeta$. It belongs to the center of $\Cg$ and 
acts on $N_k$ by conjugation $x \mapsto \zeta^{-1} x \zeta$. Then 
\begin{align*}
r_k (M) &= r_k (\tau_{\zeta}^{-1} M \tau_{\zeta})\\
&= r_k (\tau_{\zeta}^{-1}) \cdot 
{}^{\sigma_k (\tau_{\zeta}^{-1})} r_k (M \tau_{\zeta}) 
\\
&= {}^{\sigma_k (\tau_{\zeta}^{-1})} 
r_k (\tau_{\zeta})^{-1} \cdot 
\sideset{^{\sigma_k (\tau_{\zeta}^{-1})}\!}{}
{\Tmatrix{r_k (M) \cdot {}^{\sigma_k (M)} r_k (\tau_{\zeta})}}\\
&= {}^{\sigma_k (\tau_{\zeta}^{-1})} 
(r_k (\tau_{\zeta})^{-1} \cdot r_k (M) \cdot r_k (\tau_{\zeta}))\\
&= (\zeta I_{2g}) \cdot r_k (\tau_{\zeta})^{-1} \cdot 
r_k (M) \cdot r_k (\tau_{\zeta}) \cdot (\zeta^{-1} I_{2g})
\end{align*}
\noindent
where the fourth equality follows from the fact that $M$ acts on 
$N_k$ trivially. On the other hand, it is easily checked that 
\[r_k(\tau_{\zeta}) = \left( I_{2g} 
-\overline{\zegamma}^T(\overrightarrow{1-\overline{\gamma}}) \right)
(\zeta I_{2g}) \]
by using free differentials. Then 
\begin{align*}
 & \quad r_k (M) =\left( I_{2g} 
-\overline{\zegamma}^T(\overrightarrow{1-\overline{\gamma}}) \right)^{-1}
r_k (M) \left( I_{2g} 
-\overline{\zegamma}^T(\overrightarrow{1-\overline{\gamma}}) \right)\\
\Longrightarrow & \ \ \left( I_{2g} 
-\overline{\zegamma}^T(\overrightarrow{1-\overline{\gamma}}) \right)
r_k(M) = r_k(M)\left( I_{2g} 
-\overline{\zegamma}^T(\overrightarrow{1-\overline{\gamma}}) \right)\\
\Longrightarrow & \quad \overline{\zegamma}^T
(\overrightarrow{1-\overline{\gamma}})r_k(M) = 
r_k(M)\overline{\zegamma}^T(\overrightarrow{1-\overline{\gamma}}).
\end{align*}
\noindent
From (1), we see $\overline{\zegamma}^T
(\overrightarrow{1-\overline{\gamma}})r_k(M) =
\overline{\zegamma}^T(\overrightarrow{1-\overline{\gamma}})$. 
Comparing first columns, we have 
$\overline{\zegamma}^T(1-\gamma_1^{-1}) 
=r_k (M)\overline{\zegamma}^T(1-\gamma_1^{-1})$. 
(2) follows from this. 
\end{proof}
\begin{prop}\label{prop:magdet}
If $M \in \Cg [k]$ satisfies $\rank_{\KK_{N_k}}(I_{2g}-r_k(M))=2g-1$, 
then 
\[\overline{d}_{N_k}^\psi (I_{2g}-r_k(M))=
\deg^\psi\bigl(\Delta \bigl( I_{2g}-r_k(M); 
\overrightarrow{1-\overline{\gamma}}, 
\overline{\zegamma}^T \bigl) \bigl).\]
Otherwise $\overline{d}_{N_k}^\psi (I_{2g}-r_k(M))=\infty$.
\end{prop}
\begin{proof}
By Lemma \ref{uandw}, the rank of 
$I_{2g}-r_k(M)$ is at most $2g-1$. 
The case where $\rank_{\KK_{N_k}}(I_{2g}-r_k(M))<2g-1$ is 
clear by definition. Suppose 
that $\rank_{\KK_{N_k}}(I_{2g}-r_k(M))=2g-1$. 
The task is to show 
that $(\overrightarrow{1-\overline{\gamma}})$ and 
$\overline{\zegamma}^T$ satisfy the conditions (2) and (3) 
of Definition \ref{def:primitive} respectively. 
Suppose $(\overrightarrow{1-\overline{\gamma}})$ 
does not satisfy (2). 
Then $(\overrightarrow{1-\overline{\gamma}})$ can be 
written as $(\overrightarrow{1-\overline{\gamma}})=f \nu$, where 
$\nu \in (\KK_{N_k^\psi}[t^\pm])_{2g}$ and 
$f \in \KK_{N_k^\psi}[t^\pm]$ with $\deg^\psi (f) \ge 1$. Hence 
each entry of $(\overrightarrow{1-\overline{\gamma}})$ has degree 
greater than $0$. When $A \in GL(\Z N_k)$, 
the same holds for each entry of 
$(\overrightarrow{1-\overline{\gamma}})A$ which is non-zero. 
However, if we choose $f \in \Aut F_{2g}$ 
such that $\psi (f (\gamma_1))=0$, which does exist, 
\[(1-\gamma_1^{-1},\ldots,1-\gamma_{2g}^{-1})
\overline{\left(\frac{\partial f(\gamma_j)}
{\partial \gamma_i}\right)}_{i,j}
=(1-f(\gamma_1)^{-1},\ldots,1-f(\gamma_{2g})^{-1}),\]
a contradiction. Hence $(\overrightarrow{1-\overline{\gamma}})$ 
satisfies (2). 
By a similar argument, we can show that $\overline{\zegamma}^T$ 
satisfies (3), where we use $f \in \Aut F_{2g}$ preserving $\zeta$ and 
the chain rule for free differentials (see for instance 
\cite[Proposition 3.3]{bi}). 
\end{proof}
Note that $\overline{d}_{N_k}^\psi (I_{2g}-r_k(M))$ 
does not depend on the choice of the generating system of 
$\pi_1 \Sg$. This follows from the formulas 
in Proposition \ref{basischange} 
and Lemma \ref{lem:torsion} (3). 
In Section \ref{subsec:formula3}, we will see some examples showing 
that our invariants are non-trivial 
at least in the case of $k=2$.

\vspace{8pt}

\subsection{$N_k$-torsion and torsion-degree functions}
\label{subsec:torsion}
For a homology cylinder $M=(M,i_+,i_-) \in \Cg$, 
we put $\Sigma^+:=i_+(\Sg)$. 
Since the relative complex $C_\ast (M,\Sigma^+;\KK_{N_k (M)})$ 
obtained from any smooth triangulation of $(M,\Sigma^+)$ 
is acyclic by Lemma \ref{relative}, 
we can consider its Reidemeister torsion 
$\tau(C_\ast (M,\Sigma^+;\KK_{N_k (M)}))$. 
\begin{definition}
The $N_k$-{\it torsion} of a homology cylinder $M=(M,i_+,i_-) \in \Cg$ 
is given by 
\[\tau_{N_k}(M):={}^{i_+^{-1}}\tau(C_\ast (M,\Sigma^+;\KK_{N_k (M)}))
\in K_1 (\KK_{N_k})/(\pm N_k).\]
\end{definition}
\noindent
Recall that Reidemeister torsions are invariant 
under subdivision of 
the cell complex $(M,\Sigma^+)$ and simple homotopy equivalence. 

Now we consider $\tau_{N_k}(M)$ more closely. 
First we give a cell decomposition of 
$\partial M \cong \Sg \cup (-\Sg)$ by pasting two copies of 
that of $\Sg$ in Figure \ref{fig:generator}.
We denote by $R_{2g}^+$ the subcomplex $i_+(R_{2g})$. 
Take a triangulation of $\partial M$ by refining 
the cell decomposition, and extend it to one of $M$. 
Then 
\begin{align*}
\tau(C_\ast (M,R_{2g}^+;\KK_{N_k(M)})) &= 
\tau(C_\ast (\Sigma^+,R_{2g}^+;\KK_{N_k(M)})) \cdot 
\tau(C_\ast (M,\Sigma^+;\KK_{N_k(M)})) \\
&= \tau(C_\ast (M,\Sigma^+;\KK_{N_k(M)}))
\end{align*}
\noindent
by the multiplicativity of torsions and the fact 
that $\Sigma^+$ is simple homotopy equivalent to $R_{2g}^+$. 

Starting from a 3-simplex of $M$ facing the boundary, 
we can deform $M$ onto a 2-dimensional subcomplex $M'$ 
by a simple homotopy equivalence keeping 
the 1-skeleton of $M$ invariant. Then 
$\tau(C_\ast (M,R_{2g}^+;\KK_{N_k(M)}))= 
\tau(C_\ast (M',R_{2g}^+;\KK_{N_k(M)}))$. 
Next, we take a maximal tree $T$ 
of the 1-skeleton of $M'$ and collapse it to a point. 
This process also preserves 
the simple homotopy type of $(M',R_{2g}^+)$, 
so that $\tau(C_\ast (M',R_{2g}^+;\KK_{N_k(M)}))= 
\tau(C_\ast (M'/T,R_{2g}^+/(T \cap R_{2g}^+);\KK_{N_k(M)}))$.

Consequently, 
$\tau_{N_k}(M)={}^{i_+^{-1}}\tau(C_\ast 
(M'/T,R_{2g}^+/(T \cap R_{2g}^+);\KK_{N_k(M)}))$. 
$M'/T$ is a 2-dimensional cell complex with only one 
0-cell, and $R_{2g}^+/(T \cap R_{2g}^+)$ is a 
subcomplex 
consisting of one 0-cell and $2g$ 1-cells. The pair 
$(M'/T,R_{2g}^+/(T \cap R_{2g}^+))$ gives an admissible presentation 
\[\langle
i_- (\gamma_1),\ldots,i_- (\gamma_{2g}), 
z_1 ,\ldots, z_{l}, 
i_+ (\gamma_1),\ldots,i_+ (\gamma_{2g}) 
\mid r_1, \ldots r_{2g+l}\rangle\]
of $\pi_1 M$. For this presentation, we take the matrices 
$A, B, C \in M(\Z N_k)$ as in Section \ref{subsec:computation}. 
Then 
\[\tau_{N_k} (M)=
{}^{i_+^{-1}}\tau(C_\ast (M'/T,R_{2g}^+/(T \cap R_{2g}^+)
;\KK_{N_k(M)}))=
\sideset{^{i_+^{-1}\!\!}}{}{\Tmatrix{
\begin{matrix}A \\ B \end{matrix}}}.\]
Note that the matrix 
$\left( \begin{smallmatrix}A \\ B \end{smallmatrix} \right)$ 
is a presentation matrix of 
$H_1 (M'/T,R_{2g}^+/(T \cap R_{2g}^+);\Z N_k(M)) 
\cong H_1 (M,\Sigma^+;\Z N_k (M))$. 

Since multiplying an element of $\pm N_k$ does not contribute to the 
degree, we have 
\[d_{N_k}^\psi (\tau_{N_k}(M)) = 
d_{N_k}^\psi \left(
\sideset{^{i_+^{-1}\!\!}}{}{\Tmatrix{
\begin{matrix}A \\ B \end{matrix}}}\right)=
\dim_{\KK_{N_k^\psi}} 
H_1 (M,\Sigma^+;(i_+^{-1})^\ast\KK_{N_k^\psi} [t^\pm])\]
for each $\psi \in PH^1(\Sg)$. From this, we see that 
$d_{N_k}^\psi (\tau_{N_k}(M))$ 
can be computed from any admissible presentation of $\pi_1 M$. 
\begin{prop}\label{additive}
Let $M_1=(M_1,i_+,i_-), M_2=(M_2,j_+,j_-) \in \Cg$. Then 
\[d_{N_k}^\psi (\tau_{N_k}(M_1 \cdot M_2))= 
d_{N_k}^\psi (\tau_{N_k}(M_1))+d_{N_k}^{\psi \circ \sigma_2 (M_1)} 
(\tau_{N_k}(M_2))\]
holds for every $\psi \in PH^1 (\Sg)$. 
\end{prop}
\begin{proof}
We take an admissible presentation of $\pi_1 M_1$ and 
the matrices $A_{M_1}, B_{M_1}, C_{M_1}$ corresponding to it. 
We denote this presentation by 
\[\pi_1 M_1 \cong \langle i_- (\overrightarrow{\gamma}),
\overrightarrow{z},i_+ (\overrightarrow{\gamma}) \mid 
\overrightarrow{r} \rangle,\]
for short. Similarly, we take an admissible presentation 
\[\pi_1 M_2 \cong \langle j_- (\overrightarrow{\gamma}),
\overrightarrow{w},j_+ (\overrightarrow{\gamma}) \mid 
\overrightarrow{s} \rangle\]
of $\pi_1 M_2$ 
and the matrices $A_{M_2}, B_{M_2}, C_{M_2}$. 
Then we obtain an admissible presentation
\[\pi_1 (M_1 \cdot M_2)  \cong 
\langle 
j_- (\overrightarrow{\gamma}), \overrightarrow{w}, 
j_+ (\overrightarrow{\gamma}), i_- (\overrightarrow{\gamma}), 
\overrightarrow{z}, i_+ (\overrightarrow{\gamma}) \mid 
\overrightarrow{s}, 
j_+ (\overrightarrow{\gamma})i_- (\overrightarrow{\gamma})^{-1}, 
\overrightarrow{r} \rangle\]
of $\pi_1 (M_1 \cdot M_2)$. 
The corresponding partial matrix 
$\left(\begin{smallmatrix} 
A_{M_1 \cdot M_2} \\ B_{M_1 \cdot M_2}\end{smallmatrix}\right)$ 
at $\Z {N_k(M_1 \cdot M_2)}$ is given by 
\[\begin{pmatrix}
{}^{j}A_{M_2} & 0 & 0 \\
{}^{j}B_{M_2} & 0 & 0 \\
{}^{j}C_{M_2} & I_{2g} & 0 \\
0 & -I_{2g} & {}^{i}A_{M_1} \\
0 & 0 & {}^{i}B_{M_1} 
\end{pmatrix},\]
where $i:M_1 \to M_1 \cdot M_2$ and 
$j:M_2 \to M_1 \cdot M_2$ are the natural inclusions. 
From this, we have 
\begin{align*}
d_{N_k}^\psi (\tau_{N_k}(M_1 \cdot M_2)) &=
d_{N_k}^\psi \left(
\sideset{^{i \circ i_+^{-1}\!\!}}{}{\Tmatrix{
\begin{matrix}A_{M_1 \cdot M_2} \\ 
B_{M_1 \cdot M_2} \end{matrix}}}\right)\\
&=d_{N_k}^\psi \left(
\sideset{^{i_+^{-1}\!\!}}{}{\Tmatrix{
\begin{matrix}A_{M_1} \\ B_{M_1} \end{matrix}}}\right)+
d_{N_k}^\psi \left(
\sideset{^{i_+^{-1}i^{-1}j\!\!}}{}{\Tmatrix{
\begin{matrix}A_{M_2} \\ B_{M_2} \end{matrix}}}\right)\\
&=d_{N_k}^\psi \left(
\sideset{^{i_+^{-1}\!\!}}{}{\Tmatrix{
\begin{matrix}A_{M_1} \\ B_{M_1} \end{matrix}}}\right)+
d_{N_k}^\psi \left(
\sideset{^{\sigma_k (M_1) j_+^{-1}\!\!}}{}{\Tmatrix{
\begin{matrix}A_{M_2} \\ B_{M_2} \end{matrix}}}\right)\\
&=d_{N_k}^\psi (\tau_{N_k}(M_1))+
d_{N_k}^{\psi \circ \sigma_2(M_1)} (\tau_{N_k}(M_2)). 
\end{align*}
\noindent
This completes the proof. 
\end{proof}
\begin{remark}
Proposition \ref{additive} can be seen as a generalization 
of \cite[Proposition 1.11]{ld2}.
\end{remark}

\vspace{8pt}

\subsection{Factorization formulas}\label{subsec:factor}

\subsubsection{The $N_k$-degree for the closing 
of a homology cylinder}\label{subsec:formula1}
For each homology cylinder $M=(M,i_+,i_-)$, we can construct 
a closed 3-manifold defined by
\[C_M:= M/(i_+(x)=i_-(x)), \quad x \in \Sg.\]
We call it the {\it closing} of $M$. It is easily seen that 
if $M \in \Cg [k]$, we have the natural isomorphisms 
$N_k=N_k (\Sg) \cong N_k (M) \cong N_k (C_M)$. Here we identify 
these groups.

\begin{thm}\label{formula1} Let $M=(M,i_+,i_-) \in \Cg [k]$. 
For each $\psi \in PH^1 (N_k)$, we have 
\[\overline{\delta}_{N_k}^{\psi} (C_M)= 
d_{N_k}^{\psi} (\tau_{N_k} (M)) + 
\overline{d}_{N_k}^{\psi} (I_{2g}-r_k(M)) 
\in \Z_{\ge 0} \cup \{ \infty \}.\]
\end{thm}
\noindent
{\it Proof.} \ 
Take an admissible presentation of $\pi_1 M$ as in Definition 
\ref{admissible}, and construct 
the corresponding matrices $A, B, C \in M(\Z N_k)$. 

Adding $2g$ relations $i_+ (\gamma_j)=i_- (\gamma_j) \ 
(j=1,\ldots,2g)$ and deleting the generators $i_+ (\gamma_j)$ 
by using them, we obtain a presentation of $\pi_1 C_M$. 
From this presentation, we have a presentation matrix $J_{C_M}$ of 
$H_1 (C_M,p;\Z N_k)$ given by 
\[J_{C_M}=\begin{pmatrix}A+C \\ B\end{pmatrix} 
=\begin{pmatrix}I_{2g}-r_k(M) & -Z \\ 0_{(l,2g)} & I_{l}\end{pmatrix}
\begin{pmatrix}A \\ B\end{pmatrix},\]
where the second equality follows from Proposition \ref{howto}. 
Since $\left(\begin{smallmatrix}A \\ B\end{smallmatrix}\right)$ 
is invertible in $\KK_{N_k}$, 
\[\rank_{\KK_{N_k}} J_{C_M} = 
\rank_{\KK_{N_k}}\begin{pmatrix}I_{2g}-r_k(M) & -Z \\
0_{(l,2g)} & I_{l}\end{pmatrix} = 
\rank_{\KK_{N_k}} (I_{2g}-r_k(M))+l \le 2g+l-1.\]
Hence to show our claim, it suffices to prove the 
case where this value is just $2g+l-1$ 
(see Definition \ref{def:TFfunc} (2)). 

By Fundamental formula of free calculus, we have
\[(\overrightarrow{1-\overline{\gamma}} \quad 
\overrightarrow{1-\overline{z}})\ 
J_{C_M}=(1-\gamma_1^{-1},\ldots,1-\gamma_{2g}^{-1},
1-z_1^{-1},\ldots,1-z_{l}^{-1}) \ J_{C_M} =0.\]
On the other hand, we have
\[J_{C_M} \begin{pmatrix}A \\ B\end{pmatrix}^{-1}
\overline{\left(\zegamma \quad 0_{(1,l)}\right)}^T
=\begin{pmatrix}I_{2g}-r_k(M) & -Z \\ 0_{(l,2g)} & I_{l}\end{pmatrix}
\overline{\left(\zegamma \quad 0_{(1,l)}\right)}^T=0\]
by Lemma \ref{uandw} (2). 
Then we can define $\Delta \left( J_{C_M}; 
\xi, \mu\right)$, where we put 
\begin{align*}
\xi&:=(1-\gamma_1^{-1},\ldots,1-\gamma_{2g}^{-1},
1-z_1^{-1},\ldots,1-z_{l}^{-1}), \\
\mu&:=\left(\textstyle{A \atop B}\right)^{-1}
\overline{\left(\zegamma \quad 0_{(1,l)}\right)}^T. 
\end{align*}
\begin{lem}\label{lem:mu}
$\mu$ belongs to $(\Z N_k)^{2g+l}$.
\end{lem}
\begin{proof}
Recall that $\left(\begin{smallmatrix}A \\ B\end{smallmatrix}\right)$ 
is a presentation matrix of $H_1 (M, \Sigma^+;\Z N_k)$, so that 
we have an exact sequence 
\[\begin{CD}
0 @>>> (\Z N_k)^{2g+l} 
@>{\left(\begin{smallmatrix}A \\ B\end{smallmatrix}\right) \cdot}>>
(\Z N_k)^{2g+l} @>>> H_1 (M, \Sigma^+;\Z N_k) 
@>>> 0,
\end{CD}\]
where the injectivity of the second map follows from the 
fact that $H_1 (M, \Sigma^+;\KK_{N_k})=0$. 
Hence to prove the lemma, it suffices to show that 
$\overline{\left(\zegamma \quad 0_{(1,l)}\right)}^T$ 
in the third term 
$(\Z N_k)^{2g+l}=C_1 (M, \Sigma^+;\Z N_k)$ 
is mapped to $0 \in H_1 (M,\Sigma^+;\Z N_k)$. 
In the exact sequence
\[\begin{CD}
0 @>>> C_1 (\Sigma^+,p;\Z N_k)
@>>> C_1 (M, p;\Z N_k) @>>> C_1 (M, \Sigma^+;\Z N_k) 
@>>> 0, 
\end{CD}\]
the cycle $\overline{\left(\zegamma \quad 0_{(1,l)}\right)}^T$ 
is attained by 
$\overline{\left(\zegamma \quad 0_{(1,l)} \quad 0_{(1,2g)}\right)}^T 
\in C_1 (M, p;\Z N_k)=(\Z N_k)^{2g+l}$ $\oplus (\Z N_k)^{2g}$. 
Then by observing the boundary corresponding 
to the relation 
\[\prod_{j=1}^g [i_+ (\gamma_j),i_+(\gamma_{g+j})] 
\left(\prod_{j=1}^g [i_- (\gamma_j),i_-(\gamma_{g+j})]\right)^{-1},\] 
we see that 
$\overline{\left(\zegamma \quad 0_{(1,l)} \quad 0_{(1,2g)}\right)}^T$ 
is homologous to $\overline{\left(0_{(1,2g)} \quad 0_{(1,l)} 
\quad \zegamma \right)}^T$, which comes from 
$C_1 (\Sigma^+,p;\Z N_k)$. 
Our claim follows from this. 
\end{proof}

Now we continue the proof of Theorem \ref{formula1}. 
We can show that $(\xi,\mu)$ is $\psi$-primitive 
for $J_{C_M}$ as in the proof of Proposition \ref{prop:magdet}. 
Then We have 
\begin{align*}
&\Delta \left( J_{C_M}; \xi,\mu\right)=\Delta \left( J_{C_M}; 
(\overrightarrow{1-\overline{\gamma}} \quad 
\overrightarrow{1-\overline{z}}),
\left(\textstyle{A \atop B}\right)^{-1} 
\overline{\left(\zegamma \quad 0_{(1,l)}\right)}^T\right) \\
&=\Delta \left( J_{C_M} \!
\left(\textstyle{A \atop B}\right)^{-1}; 
(\overrightarrow{1-\overline{\gamma}} \quad 
\overrightarrow{1-\overline{z}}),
\overline{\left(\zegamma \quad 0_{(1,l)}\right)}^T\right) 
\cdot \det \begin{pmatrix}A \\ B\end{pmatrix}\\
&= \Delta \left( 
\left(\begin{array}{cc}
I_{2g}-r_k(M) & -Z \\ 0_{(l,2g)} & I_{l}
\end{array}\right) ; 
(\overrightarrow{1-\overline{\gamma}} \quad 
\overrightarrow{1-\overline{z}}),
\overline{\left(\zegamma \quad 0_{(1,l)}\right)}^T\right) 
\cdot \det \begin{pmatrix}A \\ B\end{pmatrix}\\
&= \Delta \left( I_{2g}-r_k(M) ; \overrightarrow{1-\overline{\gamma}}, 
\overline{\zegamma}^T \right) 
\cdot \det \begin{pmatrix}A \\ B\end{pmatrix}.
\end{align*}
From the above argument, we obtain 
\begin{align*}
\overline{\delta}_{N_k}^{\psi} (C_M) &= 
\overline{d}_{N_k}^\psi \left(J_{C_M}\right)
=\deg^\psi \left(\Delta \left( J_{C_M}; \xi,\mu\right)\right)\\
&= \deg^\psi\left(\Delta \left( I_{2g}-r_k(M); 
\overrightarrow{1-\overline{\gamma}}, 
\overline{\zegamma}^T \right) \right) 
+\deg^\psi\left( \det \begin{pmatrix}A \\ B\end{pmatrix}\right)\\
&=\overline{d}_{N_k}^{\psi} (I_{2g}-r_k(M)) 
+ d_{N_k}^{\psi} (\tau_{N_k} (M)).
\end{align*} 
\noindent
This completes the proof. \qed
\begin{remark}
When $M \in \Cg [k] \cap \Mg$, $I_{2g}-r_k(M)$ itself gives a 
presentation matrix of $H_1 (C_M,p;\Z N_k)$. Hence we have 
$\overline{\delta}_{N_k}^{\psi} (C_M)=
\overline{d}_{N_k}^{\psi} (I_{2g}-r_k(M))$, 
and moreover 
$\delta_{N_k}^{\psi} (C_M)=d_{N_k}^{\psi} (I_{2g}-r_k(M))$ 
for this case.
\end{remark}

\subsubsection{The $N_{k,T}$-degree for 
the mapping torus of a homology cylinder}\label{subsec:formula2}
Given a homology cylinder $M=(M,i_+,i_-)$, 
we have another method for obtaining a closed 3-manifold 
$T_M$ as follows. First we attach a 2-handle $I \times D^2$ 
along $I \times i_{\pm} (\partial \Sg)$, 
so that we obtain a homology 
cylinder $(M',i'_+,i_-')$ over a closed surface $\Sigma_g$, 
which corresponds to the embedding 
$\Sg \hookrightarrow \Sigma_g$. Then we put
\[T_M:=M'/(i'_+(x)=i'_-(x)), \quad x \in \Sigma_g\] 
and call $T_M$ the {\it mapping torus} of $M$. Indeed, for 
$M_{\varphi} \in \Mg \subset \Cg$, the resulting manifold 
$T_{M_{\varphi}}$ is the usual mapping torus 
of $\varphi$ extended naturally to the mapping class of $\Sigma_g$. If 
we take an admissible presentation of $\pi_1 M$ briefly denoted by 
$\langle i_- (\overrightarrow{\gamma}), 
\overrightarrow{z},i_+ (\overrightarrow{\gamma}) \mid 
\overrightarrow{r} \rangle$, 
then a presentation of $\pi_1 T_M$ is given by 
\[\pi_1 T_M \cong \langle i_- (\overrightarrow{\gamma}),
\overrightarrow{z},\lambda,i_+ (\overrightarrow{\gamma}) \mid 
\overrightarrow{r}, \textstyle\prod_{j=1}^g 
[i_-(\gamma_j),i_-(\gamma_{g+j})], 
i_- (\overrightarrow{\gamma}) \lambda 
i_+(\overrightarrow{\gamma})^{-1} \lambda^{-1} \rangle,\]
where $\lambda$ is the loop 
$I \times \{ p \} \subset I \times D^2 \subset T_M$. 
If $M \in \Cg [k]$, we have natural isomorphisms 
$N_k (\Sigma_g) \cong N_k (M')$ and 
$N_k (T_M) \cong N_k (\Sigma_g) \times \langle \lambda \rangle$. Note 
that these groups are torsion-free nilpotent. 
We consider $N_k (\Sigma_g)$ 
to be a subgroup of $N_k (T_M)$. 
For simplicity, we denote $N_k (\Sigma_g)$ by $N_{k,0}$ and 
$N_k (T_M)$ by $N_{k,T}$. 

We can show that $H_\ast (M,i_+(\Sg);\KK_{N_{k,T}})=0$ 
(see Remark \ref{acyclic}). Hence by a similar argument, 
the Magnus representation $r_{k,T}:\Cg \to GL(2g,\KK_{N_{k,T}})$ 
and the $N_{k,T}$-{\it torsion} 
\[\tau_{N_{k,T}}(M):=\tau(C_\ast (M,\Sigma^+;\KK_{N_{k,T}}))
\in K_1 (\KK_{N_{k,T}})/(\pm N_{k,T})\]
are defined. 
Then we obtain the following factorization formula of 
the $N_{k,T}$-degree 
for the mapping torus of a homology cylinder. 
\begin{thm}\label{formula2} Let $M \in \Cg [k]$. 
For each primitive element $\psi \in H^1 (N_{k,T}) 
= H^1 (T_M)$, 
the $N_{k,T}$-degree 
$\overline{\delta}_{N_{k,T}}^{\psi} (T_M)$ 
is finite, and 
we have 
\begin{align*}\overline{\delta}_{N_{k,T}}^{\psi} (T_M)&= 
\delta_{N_{k,T}}^{\psi} (T_M)\\
&=d_{N_{k,T}}^{\psi} (\tau_{N_{k,T}}(M)) 
+ d_{N_{k,T}}^{\psi} (I_{2g}-\lambda r_{k,T}(M))
-2 |\psi (\lambda)|.
\end{align*}
\end{thm}
\noindent
\begin{proof}
The first assertion is a slight generalization of 
\cite[Proposition 8.4]{har}, and we now follow the proof. 
Let $\psi \in H^1 (T_M)$ be the Poincar\'e dual of the surface 
$i'_+(\Sigma_g)=i'_-(\Sigma_g)$. This gives an exact sequence 
$1 \to N_{k,0} \to N_{k,T} \xrightarrow{\psi} \Z \to 1$.
Then our claim is proved by showing that 
$\overline{\delta}_{N_{k,T}}^{\psi} (T_M)$ is finite for 
this $\psi$. 

Let $(T_M)_{N_{k,T}}$ be the $N_{k,T}$-cover of $T_M$, and let 
$(T_M)_\psi$ be the $\Z$-cover of $T_M$ 
with respect to $\psi$. $(T_M)_\psi$ is the product 
$\cdots \cdot M' \cdot M' \cdot M' \cdot \cdots$ of 
countably many copies of $M'$, and 
$(T_M)_{N_{k,T}}$ can be 
regarded as the $N_{k,0}$-cover of $(T_M)_\psi$. Then 
\begin{align*}
H_\ast (T_M;\KK_{N_{k,T}^\psi} [t^\pm]) &=
H_\ast (C_\ast ((T_M)_{N_{k,T}}) \otimes_{N_{k,T}} \Z N_{k,T} 
(\Z N_{k,0} - \{ 0 \})^{-1})\\
&\cong H_\ast (C_\ast ((T_M)_{N_{k,T}}) \otimes_{N_{k,0}} 
\Z N_{k,0} (\Z N_{k,0} - \{ 0 \})^{-1})\\
&=H_\ast (C_\ast (((T_M)_\psi)_{N_{k,0}}) \otimes_{N_{k,0}} 
\KK_{N_{k,T}^\psi})\\
&=H_\ast ((T_M)_\psi;\KK_{N_{k,T}^\psi}).
\end{align*}
Here we remark that the image of the composite 
$\pi_1 ((T_M)_\psi) \to \pi_1 T_M \to N_{k,T}$ is contained 
in $N_{k,0}$. 
The same holds for the composite 
$\pi_1 M' \to \pi_1 T_M \to N_{k,T}$. We also remark that 
$\KK_{N_{k,T}^\psi}=\KK_{N_{k,0}}$. 

We denote by $\Sigma$ again for a lift of 
$\Sigma \subset T_M$ on $(T_M)_\psi$. 
We divide $(T_M)_\psi$ at $\Sigma$, and obtain two parts 
$(T_M)_\psi^+$ and $(T_M)_\psi^-$. 
Then $(T_M)_\psi^\pm = \varinjlim_l (M')^l$, and the inclusion 
$\Sigma \hookrightarrow (M')^l$ induces an isomorphism on 
homology. We can show that 
$H_\ast ((M')^l,\Sigma;\KK_{N_{k,T}^\psi})=0$ by the same way 
as mentioned in Lemma \ref{relative}. Thus 
$H_\ast ((T_M)_\psi^\pm, \Sigma;\KK_{N_{k,T}^\psi})=
\varinjlim_l H_\ast ((M')^l, \Sigma;$ $\KK_{N_{k,T}^\psi})=0$, 
and therefore $H_\ast ((T_M)_\psi, \Sigma;\KK_{N_{k,T}^\psi})=0$. 
This shows that 
\[H_\ast (T_M;\KK_{N_{k,T}^\psi} [t^\pm]) 
\cong H_\ast ((T_M)_\psi;\KK_{N_{k,T}^\psi}) \cong 
H_\ast (\Sigma;\KK_{N_{k,T}^\psi})\] 
is a finite dimensional 
$\KK_{N_{k,T}^\psi}$-vector space, so that 
$\overline{\delta}_{N_{k,T}}^{\psi} (T_M)$ is finite. 

To show the second assertion, we take an admissible presentation 
of $\pi_1 M$, and construct the matrices $A, B, C \in \Z N_{k,T}$ 
as before. 
From the presentation, we have a presentation matrix 
$J_{T_M}$ of $H_1 (T_M,p;\Z N_{k,T})$ given by 
\begin{align*}
J_{T_M} &=
\begin{pmatrix}
A & \overline{\zegamma}^T & I_{2g} \\
B & 0_{(l,1)} & 0_{(l,2g)}\\
0_{(1,2g+l)} & 0 & -(\overrightarrow{1-\overline{\gamma}})\\
C & 0_{(2g,1)} & -\lambda^{-1} I_{2g}
\end{pmatrix},\\
\intertext{where $\overrightarrow{\gamma}:=i_+(\overrightarrow{\gamma}) 
= i_-(\overrightarrow{\gamma})$. 
We remark that $\lambda$ belongs to the center in $N_{k,T}$. 
As presentation matrices of $H_1 (T_M,p;\Z N_{k,T})$, this matrix 
is equivalent to the square matrix}
J'_{T_M}&=\begin{pmatrix}
A+\lambda C & \overline{\zegamma}^T \\
B & 0_{(l,1)} \\
-(\overrightarrow{1-\overline{\gamma}})\lambda C & 0 
\end{pmatrix}.\\
\intertext{By Proposition \ref{howto}, Lemma \ref{uandw} and 
Fundamental formula of free calculus, we have} 
J'_{T_M}&=\begin{pmatrix}
A+\lambda C & \overline{\zegamma}^T \\
B & 0_{(l,1)} \\
\lambda (\overrightarrow{1-\overline{\gamma}} \quad 
\overrightarrow{1-\overline{z}})
\left(\textstyle{A \atop B}\right)& 0 
\end{pmatrix}\\
&=\begin{pmatrix}
I_{2g}-\lambda r_{k,T} (M) & - \lambda Z & \overline{\zegamma}^T\\
0_{(l,2g)} & I_l & 0_{(l,1)}\\
\lambda (\overrightarrow{1-\overline{\gamma}}) & 
\lambda (\overrightarrow{1-\overline{z}}) & 0
\end{pmatrix}
\begin{pmatrix}
A & 0_{(2g,1)}\\
B & 0_{(l,1)} \\ 
0_{(1,2g+l)} & 1 
\end{pmatrix}.
\end{align*}
Note that $\widetilde{\left(A \atop B\right)}:=
\left(\begin{smallmatrix}
A & 0_{(2g,1)}\\
B & 0_{(l,1)} \\ 
0_{(1,2g+l)} & 1 
\end{smallmatrix}\right)$ is invertible in $\KK_{N_{k,T}}$. 
Then it is easily checked that 
\begin{align*}
&(\overrightarrow{1-\overline{\gamma}} \quad 
\overrightarrow{1-\overline{z}} \quad 1-\lambda^{-1})\ 
J'_{T_M}=0,\\
&J'_{T_M} \widetilde{\textstyle\left(A \atop B\right)}^{-1} 
\overline{(\zegamma \quad 0_{(1,l)} \quad \lambda^{-1}-1)}^T =0.
\end{align*}
We put $\widetilde{\xi}:=(\overrightarrow{1-\overline{\gamma}} \quad 
\overrightarrow{1-\overline{z}} \quad 1-\lambda^{-1})$ and 
$\widetilde{\mu}:=\widetilde{\textstyle\left(A \atop B\right)}^{-1} 
\overline{(\zegamma \quad 0_{(1,l)} \quad \lambda^{-1}-1)}^T$. 
As in Proposition \ref{prop:magdet} and Lemma \ref{lem:mu}, 
we can show that $(\widetilde{\xi},\widetilde{\mu})$ is 
$\psi$-primitive for $J'_{T_M}$. Then
\begin{align*}
\delta_{N_{k,T}}^{\psi}  (T_M) =& d_{N_{k,T}}^\psi 
\left( J'_{T_M}\right)=
\deg^\psi \left( \Delta \left( J'_{T_M};
\widetilde{\xi},\widetilde{\mu}\right) \right)\\
=& \deg^\psi \left( \Delta \left( J'_{T_M}
\widetilde{\textstyle\left(A \atop B\right)}^{-1};\widetilde{\xi},
\widetilde{\textstyle\left(A \atop B\right)}\widetilde{\mu}\right) 
\cdot \det \widetilde{\textstyle\left(A \atop B\right)} \right)\\
=& \deg^\psi \left( \Delta \left( 
\begin{pmatrix}
I_{2g}-\lambda r_{k,T} (M) & - \lambda Z & \overline{\zegamma}^T\\
0_{(l,2g)} & I_l & 0_{(l,1)}\\
\lambda (\overrightarrow{1-\overline{\gamma}}) & 
\lambda (\overrightarrow{1-\overline{z}}) & 0
\end{pmatrix}; \widetilde{\xi}, 
\widetilde{\textstyle\left(A \atop B\right)}\widetilde{\mu}\right) 
\cdot \det \widetilde{\textstyle\left(A \atop B\right)} \right)\\
=& \deg^\psi \left( \det 
\begin{pmatrix} 
I_{2g}-\lambda r_{k,T} (M) & - \lambda Z\\
0_{(l,2g)} & I_l
\end{pmatrix}\cdot \lambda^{-1}(1-\lambda^{-1})^{-2} \cdot 
\det \widetilde{\textstyle\left(A \atop B\right)} \right)\\
=& \deg^\psi \left( \det(I_{2g}-\lambda r_{k,T} (M)) 
\cdot \lambda^{-1}(1-\lambda^{-1})^{-2} \cdot 
\det \widetilde{\textstyle\left(A \atop B\right)} \right)\\
=& \deg^\psi (\det( I_{2g}-\lambda r_{k,T}(M))) + 
\deg^\psi (\det(\tau_{N_{k,T}} (M)))-2 |\psi (\lambda)|.
\end{align*}
This completes the proof. 
\end{proof}

\subsubsection{The case of $k=2$ $($commutative case$)$}
\label{subsec:formula3}
Since $\Z N_2=\Z N_2 (\Sigma_g)$ and $\KK_{N_2}=\KK_{N_2 (\Sigma_g)}$ 
are commutative, 
we can use the ordinary determinant for computation. 
Moreover, we can obtain some invariants before taking degrees. 
For example, define
\[\Delta(M):=(-1)^{i+j}
\frac{\det\left( (I_{2g}-r_2 (M))_{(i,j)}\right)}
{(1-\gamma_i^{-1})(\overline{\frac{\partial \zeta}{\partial \gamma_j}})} 
\in \KK_{N_2},\]
where $A_{(i,j)}$ is the matrix obtained from a matrix $A$ by 
removing its $i$-th row and $j$-th column. 
$\Delta(M)$ is well-defined by Lemma \ref{lem:torsion}. 
Note that this invariant is based on that for 
string links given in \cite{klw}, and 
we call it the {\it Alexander rational function} of $M$. 

\begin{thm}\label{decomp}
Let $M \in \Cg [2]$, and let $\Delta_{C_M}$, 
$\Delta_{T_M}$ be the Alexander polynomials of 
$C_M$, $T_M$, respectively. Then 
\begin{align*}
\Delta_{C_M} &\doteq \overline{\tau_{N_2}(M) \cdot \Delta(M)}, \\
\Delta_{T_M} &\doteq \overline{\tau_{N_2}(M) \cdot 
\det \bigl(I_{2g}-\lambda r_{2,T} (M) \bigl) 
\cdot (1-\lambda^{-1})^{-2}},
\end{align*}
where $\doteq$ means that these equalities hold in 
$\KK_{N_2}$ and $\KK_{N_2 (T_M)}$ up to $\pm N_2$ and $\pm N_2 (T_M)$ 
respectively. 
\end{thm}
\noindent
\begin{proof}
We prove only the first assertion. The proof is almost the same as 
that for Theorem \ref{formula1} under the following remarks. 
We follow the notation used there. We may assume that 
$\rank_{\KK_{N_2}} J_{C_M} = \rank_{\KK_{N_2}} 
r_2 (M)+l=2g+l-1$. 

By definition, 
$\Delta_{C_M}$ is the greatest common divisor of 
$\{ \det \overline{J_{C_M \, (i,j)}}^T \}_{1 \le i,j \le 2g+l}$. 
We show that it is nothing other than 
\[\Delta:=\overline{\Delta \left( J_{C_M}; 
(\overrightarrow{1-\overline{\gamma}} \quad 
\overrightarrow{1-\overline{z}}),
\left(\textstyle{A \atop B}\right)^{-1} 
\overline{\left(\zegamma \quad 0_{(1,l)}\right)}^T\right)}
\doteq \overline{\det \left(A \atop B\right) \cdot \Delta(M)}.\]
As seen in Lemma \ref{lem:mu}, 
$\left(\textstyle{A \atop B}\right)^{-1} 
\overline{\left(\zegamma \quad 0_{(1,l)}\right)}^T$ 
is a vector in $(\Z N_2)^{2g+l}$. 
If $\Delta$ is in $\Z N_2$, 
it attains the greatest common divisor. 
To show it, 
suppose $\Delta=h_1/h_2$ where $h_1 \in \Z N_2$ 
and $h_2 \in \Z N_2 -\{0\}$ 
are relatively prime. From the definition of $\Delta$, 
we have 
\[\frac{(1-\gamma_i^{-1})
\left(\overline{\frac{\partial \zeta}{\partial \gamma_j}}\right)h_1}{h_2}
=(-1)^{i+j}\det \overline{J_{C_M \, (i,j)}}^T \in \Z N_2.\] 
Hence $h_2$ is a common divisor of $\left\{(1-\gamma_i^{-1})
\left(\overline{\frac{\partial \zeta}{\partial \gamma_j}}\right)
\right\}_{i,j}$'s, 
and it is $1$. That is, $h_2$ is a unit in $\Z N_2$. 

$\det \left(A \atop B\right) \in \KK_{N_2}$ 
(up to $\pm N_2$) does not 
depend on the choice of an admissible presentation, and 
it gives $\tau_{N_2} (M)$. 
Indeed the matrix $\left(A \atop B\right)$ 
is a presentation matrix 
of $H_1 (M, \Sigma^+;\Z N_2)$, and its determinant gives 
a generator of the 0-th elementary ideal, which is principal and 
invariant under Tietze transformations. This completes the proof.
\end{proof}

The formula in Theorem \ref{formula1} holds as elements of 
$\Z \cup \{ \infty \}$, so that the additivity loses its 
meaning when the value is $\infty$. Note that 
$\overline{\delta}_{N_k}^{\psi}(C_M)=\infty$ if and only if 
$\overline{d}_{N_k}^{\psi} (I_{2g}-r_k(M))=\infty$, 
and this occurs when 
$H_1 (C_M;\KK_{N_k^\psi} [t^\pm])$ has a non-trivial free 
part. 
The following are some examples of homology cylinders 
which have non-trivial Alexander rational functions. 
By using Theorem \ref{harvey} in the next subsection, 
we obtain many situations where the formula sufficiently works. 
When $k \ge 3$, the computation becomes quite difficult in general. 

\begin{example}\label{ex:mag1}
Assume that $g=1$. The Dehn twist $\tau_{\zeta} \in 
\mathcal{M}_{1,1}$ belongs to $\mathcal{C}_{1,1} [3]$. 
Then, we have 
\[r_2(\tau_{\zeta}) = 
\begin{pmatrix}
\gamma_1^{-1} +\gamma_2^{-1} -\gamma_1^{-1} \gamma_2^{-1} & 
-1+2\gamma_2^{-1} -\gamma_2^{-2}\\
1-2\gamma_1^{-1} +\gamma_1^{-2} & 
2 -\gamma_1^{-1} -\gamma_2^{-1} +\gamma_1^{-1} \gamma_2^{-1}
\end{pmatrix}.\]
Then $\Delta(\tau_{\zeta})=1 \in \Z N_2$, 
which is non-trivial. 
\end{example}
\begin{example}\label{ex:mag2}
Assume that $g \ge 2$. Let $\tau_1$, $\tau_2$ and 
$\tau_3$ be Dehn twists along simple closed curves $c_1$, 
$c_2$ and $c_3$ as in Figure \ref{fig:eg2}.

\begin{figure}[htbp]
\begin{center}
\includegraphics{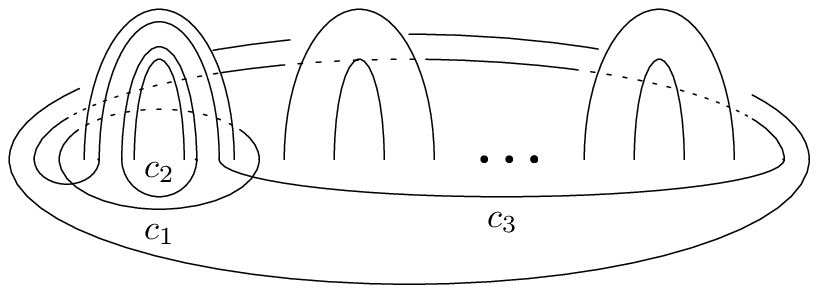}
\end{center}
\caption{}
\label{fig:eg2}
\end{figure}

\noindent
Then $\tau_1 \tau_2^{-1}, \tau_3 \in \Cg [2]$. 
By a direct computation, we can check that 
$\Delta(\tau_1 \tau_2^{-1} \cdot \tau_3) =
-(\gamma_1-1)^{2g-2}$,
although $\Delta(\tau_1 \tau_2^{-1})=\Delta(\tau_3)=0$. 
\end{example}

\subsection{$N_k$-torsions and Harvey's Realization Theorem}
\label{subsec:Harvey}
By Proposition \ref{additive}, the degree of 
the $N_k$-torsion gives a monoid homomorphism 
\[d_{N_k}^\psi (\tau_{N_k} ( \cdot )):
\Cg[2] \longrightarrow \Z_{\ge 0}\] 
for each $\psi \in PH^1 (\Sg)$ and an integer $k \ge 2$. 
To see some properties of these homomorphisms, including their 
non-triviality, we use a variant of 
Harvey's Realization Theorem in \cite[Theorem 11.2]{har} 
which gives a method for performing surgery on a 
compact orientable 3-manifold to obtain 
a homology cobordant one having distinct higher-order degrees. 
\begin{thm}\label{harvey}
Let $M \in \Cg$ be a homology cylinder. For each 
primitive element $x$ of $H_1 (\Sg)$ and 
any integers $n \ge 2$ and $k \ge 1$, there 
exists a homology cylinder $M'$ such that 
\begin{enumerate}
\item $M'$ is homology cobordant to $M$, 
\item $d_{N_l}^\psi (\tau_{N_l}(M')) = 
d_{N_l}^\psi (\tau_{N_l}(M))$ \quad 
for $2 \leq l \le n-1$, 
\item $d_{N_n}^\psi (\tau_{N_n}(M')) \ge 
d_{N_n}^\psi (\tau_{N_n}(M)) +k | p |$
\end{enumerate}
\noindent
for any $\psi \in PH^1 (\Sg)$ satisfying 
$\psi (x) =p$.
\end{thm}
\begin{proof}
The proof is based on Harvey's proof of Realization Theorem 
in \cite[Theorem 11.2]{har}. However, 
since we now use the lower central series instead of 
the rational derived series, we can shorten the argument. 

We take a loop representing $x \in H_1 (\Sg)$, 
and denote it by $x$ again. 
We also take a loop $\gamma$ whose homology class 
in $H_1 (\Sg)$ is independent of $x$. 

We attach a 1-handle to $M \times \{ 1\} \subset 
M \times I$, and then attach a 2-handle to obtain a 
4-manifold $W$. Here the 2-handle are attached along the loop 
$\alpha [X_{n-1}, A_{k+1}]$, where 
$\alpha \in \pi_1 M$ is a loop corresponding to the added 1-handle, 
and $X_{n-1}, A_{k+1} \in \pi_1 M$ are inductively defined by
\[\begin{array}{ll}
X_1 = i_+ (x), & 
X_l = [i_+ (\gamma), X_{l-1}] \ \ \mbox{for $l \ge 2$}, \\
A_1 = \alpha, & A_l = [i_+ (x), A_{l-1}] \ \ \mbox{for $l \ge 2$}. 
\end{array}\]
It is easily seen that 
$X_l \in \Gamma^l (\pi_1 M) - \Gamma^{l+1} (\pi_1 M)$. 
$M'$ is defined as another part of $\partial W$, namely 
$\partial W = M \cup M'$ and $M \cap M'=\partial M =\partial M'$. 
From the construction, we have $H_\ast (W,M)=0$. We also have 
$H_\ast (W,M')=0$ by using 
the Poincar\'e-Lefschetz duality and the 
universal coefficient theorem. Hence $(M', i_+, i_-) \in \Cg$, 
and it is homology cobordant to $M$. 
Stallings' theorem shows that 
$N_l \xrightarrow{i_+} N_l (M) \to 
N_l (W) \gets N_l (M') \xleftarrow{i_+} N_l$ are all 
isomorphisms. Using them, we identify 
$N_l, N_l (M), N_l (M')$ and $N_l (W)$. 

For simplicity, we put 
$K_l := \KK_{N_l^\psi}[t^\pm]=\Z N_l (\Z N_l^\psi - \{0\})^{-1}$. 
Recall that $H_\ast (M, \Sigma^+;\KK_{N_l}) 
=H_\ast (M', \Sigma^+;\KK_{N_l})=0$ as in Lemma \ref{relative}. 
By the same proof, we have $H_\ast (W, \Sigma^+;\KK_{N_l})=0$. 
Hence $H_\ast (M, \Sigma^+;K_l), 
H_\ast (M', \Sigma^+;K_l)$ and $H_\ast (W, \Sigma^+;K_l)$ are 
all finite dimensional $\KK_{N_l^\psi}$-vector spaces. 
As seen in Section \ref{subsec:torsion}, 
$d_{N_l}^\psi (\tau_{N_l}(M)) = \dim_{\KK_{N_l^\psi}} 
H_1 (M,\Sigma^+;K_l)$. If we take an admissible presentation 
of $\pi_1 M$ and the matrices $A, B \in \Z N_k$ as before, 
$\left(A \atop B\right)$ gives a presentation matrix 
of $H_1 (M,\Sigma^+;K_l)$. Then one of $H_1 (W, \Sigma^+;K_l)$ is 
given by
\[\left(\begin{array}{cc}
A \atop B & \ast\\ 
0_{(1,2g+l)} & \overline{\textstyle\frac{\partial \alpha 
[X_{n-1}, A_{k+1}]}{\partial \alpha}}
\end{array}\right),\]
so that 
\[\dim_{\KK_{N_l^\psi}} H_1 (W, \Sigma^+;K_l) = 
d_{N_l}^\psi (\tau_{N_l}(M))+\deg^\psi \left(
\overline{\textstyle\frac{\partial \alpha 
[X_{n-1}, A_{k+1}]}{\partial \alpha}}\right).\] 
By a direct computation, 
\[\textstyle\frac{\partial \alpha 
[X_{n-1}, A_{k+1}]}{\partial \alpha} =
1+\alpha \left\{ (1-X_{n-1}A_{k+1}X_{n-1}^{-1}) 
\textstyle\frac{\partial X_{n-1}}{\partial \alpha}+
(X_{n-1}-[X_{n-1},A_{k+1}])
\textstyle\frac{\partial A_{k+1}}{\partial \alpha}\right\}.\]
When $2 \le l \le n-1$, we have $X_{n-1}=A_{k+1}=1 \in N_l$, 
so that $\textstyle\frac{\partial \alpha 
[X_{n-1}, A_{k+1}]}{\partial \alpha}=1$, and 
$H_1 (M, \Sigma^+;K_l) \cong H_1 (W, \Sigma^+;K_l)$. 
When $l=n$, we have 
$X_{n-1} \neq A_{k+1}=1 \in N_l$, so that 
\begin{align*}
\textstyle\frac{\partial \alpha 
[X_{n-1}, A_{k+1}]}{\partial \alpha}&=1+(X_{n-1}-1)
\textstyle\frac{\partial A_{k+1}}{\partial \alpha}
=1+(X_{n-1}-1)(x-A_{k+1})
\textstyle\frac{\partial A_k}{\partial \alpha}\\
&= \cdots =1+(X_{n-1}-1)(x-A_{k+1})(x-A_k)\cdots (x-A_2), 
\end{align*}
\noindent
and 
\[\deg^\psi \left(
\overline{\textstyle\frac{\partial \alpha 
[X_{n-1}, A_{k+1}]}{\partial \alpha}}\right)=
\deg^\psi \left(\textstyle\frac{\partial \alpha 
[X_{n-1}, A_{k+1}]}{\partial \alpha}\right)=
\begin{cases}k |p| & (n \ge 3)\\ (k+1) |p| & (n=2)
\end{cases}.\]
In each case, 
$\dim_{\KK_{N_n^\psi}} H_1 (W, \Sigma^+;K_n) \ge  
d_{N_n}^\psi (\tau_{N_n}(M))+k |p|$.

By considering the dual handle decomposition, we see that 
$W$ is obtained from $M' \times I$ by attaching a 2-handle and 
a 3-handle. Hence 
$H_1 (M', \Sigma^+;K_l) \to H_1 (W, \Sigma^+;K_l)$ 
is an epimorphism. 
In particular, when $l=n$, 
\[d_{N_n}^\psi (\tau_{N_n}(M')) \ge 
\dim_{\KK_{N_n^\psi}} H_1 (W, \Sigma^+;K_n) \ge  
d_{N_n}^\psi (\tau_{N_n}(M))+k |p|.\]

It remains to proof that the map $H_1 (M', \Sigma^+;K_l) 
\to H_1 (W, \Sigma^+;K_l)$ is injective when $2 \le l \le n-1$. 
We now show that $H_2 (W,M';K_l)=0$. 
By the Poincar\'e-Lefschetz duality, 
$H_2 (W,M';K_l) \cong H^2 (W,M;K_l)$. On the other hand, 
it is easily checked that $H_0 (W,M;K_l)=H_1 (W,M;K_l)
=H_2 (W,M;K_l)=0$. Then the universal coefficient spectral sequence 
(see \cite[Theorem 2.3]{le0}) shows our claim. 
Consequently, 
$H_1 (M, \Sigma^+;K_l) \cong H_1 (W, \Sigma^+;K_l) 
\cong H_1 (M', \Sigma^+;K_l)$ and 
$d_{N_l}^\psi (\tau_{N_l}(M'))=d_{N_l}^\psi (\tau_{N_l}(M))$. 
This completes the proof.
\end{proof}
\begin{cor}
For any $\psi \in PH^1 (\Sg)$, 
the maps $\{d_{N_k}^\psi (\tau_{N_k} ( \cdot )):
\Cg[2] \to \Z_{\ge 0}\}_{k \ge 2}$ 
are all non-trivial homomorphisms, and independent of each other. 
\end{cor}
\noindent
In fact, we can show it by constructing homology cylinders 
that are homology cobordant to the unit $1_{\Cg}$. From this 
we see that $\Cg [2], \Cg[3],\ldots,\Ker(\Cg \to \Hg)$ are not 
finitely generated monoids. 
Note that $d_{N_k}^\psi (\tau_{N_k} (M))=0$ if 
$M \in \Mg$, since $\Sg \times I$ is simple homotopy equivalent to 
$\Sg$ and hence $\tau_{N_k} (M)$ is trivial. 

\subsection{Appendix:Application of torsion-degree functions 
to $\Aut \Acy_n$}\label{subsec:acyclic}
In \cite{sa2}, we defined the Magnus representation 
$r_k : \Aut \Acy_n \to GL(n,\KK_{N_k (F_n)})$ for 
$\Aut \Acy_n$, where $\Acy_n$ is a completion of $F_n$ in 
a certain sense and is called the {\it acyclic closure} of $F_n$. 
The natural map $F_n \to \Acy_n$ is known to be injective and 2-connected. 
In particular, $N_k (F_n)= N_k (\Acy_n)$, and 
we denote it briefly by $N_k$ in this subsection. 
$\Aut \Acy_n$ can be regarded as an enlargement of $\Aut F_n$. Indeed 
we have the enlarged Dehn-Nielsen homomorphism 
$\sigma^{\mathrm{acy}}:\Hg \to \Aut \Acy_{2g}$ extending the 
classical one $\sigma:\Mg \hookrightarrow \Aut F_{2g}$. 
That is, we have the commutative diagram
\[\begin{array}{ccc}
\Aut F_{2g} & \hookrightarrow & \Aut \Acy_{2g}\\
\mbox{\rotatebox[origin=c]{90}{$\hookrightarrow$}{\tiny$\sigma$}}
&&\mbox{$\uparrow${\tiny$\sigma^{\mathrm{acy}}$}}\\
\Mg & \hookrightarrow & \Hg
\end{array}.\]
Note that $\sigma^{\mathrm{acy}}$ is not injective. 
The Magnus representation for homology cylinders 
is nothing other than the composite 
$\Hg \stackrel{\sigma^{\mathrm{acy}}}{\longrightarrow} \Aut \Acy_{2g} 
\xrightarrow{r_k} GL(2g,\KK_{N_k})$. 

We now consider the map $d_{N_k}^\psi \circ r_k: \Aut \Acy_n 
\to \Z$ for $\psi \in PH^1(F_n)$, where $PH^1(F_n)$ denotes the set of 
primitive elements of $H^1(F_n)$. 
Since $d_{N_k}^\psi (A)=0$ for $A \in GL(\Z N_k)$, 
it follows that $d_{N_k}^\psi \circ r_k \big|_{\Aut F_n}$ is trivial. 
When $n=2g$, 
$d_{N_k}^\psi \circ r_k \big|_{\Im \sigma^{\mathrm{acy}}}$ is also 
trivial as seen in Theorem \ref{vanish}. On the other hand, 
$d_{N_k}^\psi \circ r_k$ is actually non-trivial 
on $\Aut \Acy_n$ as we will see below. Since $r_k$ is a crossed 
homomorphism, we have the following.
\begin{prop}\label{prop:autacy1}
For $f, g \in \Aut \Acy_n$ and $\psi \in PH^1 (F_n)$, we have
\[d_{N_k}^\psi (r_k (fg))=d_{N_k}^\psi (r_k (f))+
d_{N_k}^{\psi \circ f} (r_k (g)).\]
\end{prop}
\noindent
In particular, if we restrict $d_{N_k}^\psi \circ r_k$ 
to $\IAut \Acy_n := 
\Ker (\Aut \Acy_n \to \Aut N_2=GL(n,\Z))$, it becomes a homomorphism. 
\begin{remark}
$\Aut \Acy_n$ acts on $PH^1(F_n)$ from the right, and hence 
acts on $\mathrm{Map}(PH^1(F_n),\Z)$ from the left. 
We regard $d_{N_k}^{\,\cdot} (r_k (\cdot))$ as 
a map $\Aut \Acy_n \to \mathrm{Map}(PH^1(F_n),\Z)$. Then 
Proposition \ref{prop:autacy1} shows that 
$d_{N_k}^{\,\cdot} (r_k (\cdot))$ is a 1-cocycle in 
$C^1 (\Aut \Acy_n, \mathrm{Map}(PH^1(F_n),\Z))$. 
We can see that it is non-trivial in 
$H^1 (\Aut \Acy_n, \mathrm{Map}(PH^1(F_n),\Z))$ from the proof of 
Theorem \ref{thm:autacy2} below. 
\end{remark}
\begin{thm}\label{thm:autacy2}
For every $n \ge 2$, $\IAut \Acy_n$ is not finitely generated. 
In fact, $H_1(\IAut \Acy_n)$ has infinite rank.
\end{thm}
\begin{proof}
Let $F_n =\langle x_1,x_2, \ldots, x_n \rangle$. 
We take $\psi:=x_1^\ast \in PH^1 (F_n)$. Consider the 
endomorphism $f_k$ of $F_n$ given by 
\[f_k(x_1)=x_1[Y_{k-1},Y_k], \quad f_k(x_i)=x_i 
\ \ \mbox{for $i \ge 2$},\]
where we define $Y_1 =x_1$ and $Y_l=[x_2, Y_{l-1}]$ for $l \ge 2$. 
Since $f_k$ is 2-connected, it induces an 
automorphism of $\Acy_n$ (see \cite[Section 4]{sa2}). 
We denote it by $f_k$ again. It belongs to $\IAut \Acy_n$. 
For such an automorphism, the Magnus matrix $r_l (f_k)$ can be computed 
by using free differentials. That is, we have 
\[r_l(f_k)=\begin{pmatrix}
\overline{\textstyle\frac{\partial f(x_1)}{\partial x_1}} & 0_{(1,n-1)}\\
\begin{smallmatrix}
\overline{\textstyle\frac{\partial f(x_1)}{\partial x_2}}\\
\vdots \\ 
\overline{\textstyle\frac{\partial f(x_1)}{\partial x_n}}
\end{smallmatrix} & I_{n-1}
\end{pmatrix}\]
at $\Z N_k$. Then 
$d_{N_l}^\psi(r_l(f_k))=\deg^\psi(\det(r_l(f_k)))=\deg^\psi 
\left(\overline{\textstyle\frac{\partial f(x_1)}{\partial x_1}}
\right)=\deg^\psi 
\left(\textstyle\frac{\partial f(x_1)}{\partial x_1}\right)$. 
By a direct computation, we have 
\[\textstyle\frac{\partial f(x_1)}{\partial x_1} =
1+x_1 \left\{ (1-Y_{k-1}Y_k Y_{k-1}^{-1}) 
\textstyle\frac{\partial Y_{k-1}}{\partial x_1}+
(Y_{k-1}-[Y_{k-1},Y_k])\textstyle\frac{\partial Y_k}
{\partial x_1}\right\}.\]
When $2 \le l \le k-1$, we have $Y_{k-1}=Y_k=1 \in N_l$, 
so that 
\[d_{N_l}^\psi (r_l(f_k))=\deg^\psi\left(\textstyle\frac{\partial f(x_1)}
{\partial x_1}\right)=\deg^\psi(1)=0.\] 
When $l=k$, we have $Y_{k-1} \neq 1 \in N_k$, so that 
\begin{align*}
\textstyle\frac{\partial f(x_1)}{\partial x_1}&=
1+x_1(Y_{k-1}-1)
\textstyle\frac{\partial Y_k}{\partial x_1}
=1+x_1(Y_{k-1}-1)(x_2-Y_k)
\textstyle\frac{\partial Y_{k-1}}{\partial x_1}\\
&= \cdots =1+x_1(Y_{k-1}-1)(x_2-Y_k)(x_2-Y_{k-1})\cdots (x_2-Y_2), 
\end{align*}
\noindent
and 
\[d_{N_k}^\psi (r_k(f_k))=\deg^\psi 
\left(\textstyle\frac{\partial f(x_1)}{\partial x_1}\right)=
\begin{cases}1 & (k \ge 3)\\ 2 & (k=2)
\end{cases}.\]
This shows that $\{d_{N_k}^\psi (r_k(\cdot))\}_{k \ge 2}$ 
are all non-trivial, and independent of each other. Our claim 
follows from this.
\end{proof}

\vspace{10pt}

\section{Acknowledgement}\label{sec:acknowledge}
The author would like to express his gratitude 
to Professor Shigeyuki Morita 
for his encouragement and helpful suggestions. 
He also would like to thank Professor Masaaki Suzuki for valuable 
discussions and advice.

This research was  supported by 
JSPS Research Fellowships for Young Scientists.


\begin{thebibliography}{aaaa}
\bibitem{bi} J.~Birman, 
\textit{Braids, Links and Mapping Class Groups}, Ann.\ of Math.\ 
Stud. 82, Princeton Univ.\ Press (1974)

\bibitem{coc} T.~Cochran, 
\textit{Noncommutative knot theory}, 
Algebr.\ Geom.\ Topol. 4 (2004), 347--398

\bibitem{co2} P.~M.~Cohn, 
\textit{Skew Fields; Theory of general division rings}, 
Encyclopedia Math.\ Appl. Cambridge Univ.\ Press, Cambridge 
(1995)

\bibitem{fo} R.~H.~Fox, 
\textit{Free differential calculus. II}, 
Ann.\ of Math. 59(2) (1954), 196--210

\bibitem{fr} S.~Friedl, 
\textit{Reidemeister torsion, the Thurston norm 
and Harvey's invariants}, preprint

\bibitem{gl} S.~Garoufalidis, J.~Levine, 
\textit{Tree-level invariants of three-manifolds, Massey products 
and the Johnson homomorphism}, 
Graphs and patterns in mathematics and 
theoretical physics, 
Proc.\ Sympos.\ Pure Math. 73 (2005), 173--205 

\bibitem{ha} K.~Habiro, 
\textit{Claspers and finite type invariants of links}, 
Geom.\ Topol. 4 (2000), 1--83

\bibitem{har} S.~Harvey, 
\textit{Higher-order polynomial invariants of 3-manifolds 
giving lower bounds for the Thurston norm}, 
Topology 44 (2005), 895--945

\bibitem{har2} S.~Harvey, 
\textit{Monotonicity of degrees of generalized 
Alexander polynomials of groups and 3-manifolds}, 
Math.\ Proc.\ Cambridge Philos.\ Soc. 140 (2006), 431--450

\bibitem{klw} P.~Kirk, C.~Livingston, Z.~Wang, 
\textit{The Gassner representation for string links}, 
Commun.\ Contemp.\ Math. 1(3) (2001), 87--136

\bibitem{ld} J.~Y.~Le\ Dimet, 
\textit{Enlacements d'intervalles et repr\'esentation de Gassner}, 
Comment.\ Math.\ Helv. 67 (1992), 306--315

\bibitem{ld2} J.~Y.~Le\ Dimet, 
\textit{Enlacements d'intervalles et torsion de Whitehead}, 
Bull.\ Soc.\ Math.\ France 129 (2001), 215--235

\bibitem{le0} J.~Levine, 
\textit{Knot modules I}, 
Trans.\ Amer.\ Math.\ Soc. 229 (1977), 1--50

\bibitem{le} J.~Levine, 
\textit{Homology cylinders: an enlargement of the mapping 
class group}, 
Algebr.\ Geom.\ Topol. 1 (2001), 243--270

\bibitem{mi} J.~Milnor, 
\textit{Whitehead torsion}, 
Bull.\ Amer.\ Math.\ Soc 72 (1966), 358--426

\bibitem{mo} S.~Morita, 
\textit{Abelian quotients of subgroups of the mapping class 
group of surfaces}, Duke Math.\ J. 70 (1993), 699--726

\bibitem{pa} D.~Passman, 
\textit{The Algebraic Structure of Group Rings}, 
John Wiley and Sons (1975)

\bibitem{sa2} T.~Sakasai, 
\textit{Homology cylinders and the acyclic closure of a free group}, 
Algebraic \& Geometric Topology 6 (2006), 603--631

\bibitem{sa5} T.~Sakasai, 
\textit{The symplecticness of the Magnus representation 
for homology cobordisms of surfaces}, in preparation

\bibitem{st} J.~Stallings, 
\textit{Homology and central series of groups}, 
J.\ Algebra 2 (1965), 170--181

\bibitem{tu2} V.~Turaev, 
\textit{Introduction to combinatorial torsions}, Lectures 
Math.\ ETH Z\"urich, Birkh\"auser (2001)

\end{thebibliography}
\end{document}